\documentclass[a4paper]{article}

\usepackage[a4paper,left=3cm,right=3cm,top=3cm,bottom=3cm]{geometry}
\usepackage{amsmath}
\usepackage{amssymb}
\usepackage[usenames]{color}
\usepackage{booktabs}
\usepackage{multirow}
\usepackage{comment}
\usepackage{tabularx}
\usepackage{graphicx} 
\usepackage{subcaption}
\usepackage[ruled,vlined]{algorithm2e}
\usepackage{hyperref}
\usepackage{enumerate}
\usepackage{authblk}

\newcommand{\cA}{\mathcal{A}}
\newcommand{\cD}{\mathcal{D}}
\newcommand{\cL}{\mathcal{L}}
\newcommand{\R}{\mathbb{R}}

\DeclareMathOperator*{\argmin}{argmin}

\usepackage[style=numeric-comp,sorting=nyt,backend=biber,giveninits=true]{biblatex}
\title{Algorithmic unfolding for image reconstruction and localization problems in fluorescence microscopy}
\author[1]{Silvia Bonettini\thanks{silvia.bonettini@unimore.it}}
\author[2]{Luca Calatroni\thanks{calatroni@i3s.unice.fr}}
\author[1]{Danilo Pezzi\thanks{danilo.pezzi@unimore.it}}
\author[1]{Marco Prato\thanks{marco.prato@unimore.it}}
\affil[1]{Dipartimento di Scienze Fisiche, Informatiche e Matematiche, Università di Modena e Reggio Emilia, Via Campi 213/b, 41125, Modena, Italy}
\affil[2]{CNRS, Laboratoire I3S, UMR 7271, UniCa, 2000,  Route des Lucioles, Sophia-Antipolis, 06903, France}

\date{}

\begin{document}

\maketitle

\bigskip

\begin{abstract}
We propose an unfolded accelerated projected-gradient descent procedure to estimate model and algorithmic parameters for image super-resolution and molecule localization problems in image microscopy. The variational lower-level constraint enforces sparsity of the solution and encodes different noise statistics (Gaussian, Poisson), while the upper-level cost assesses optimality w.r.t.~the task considered. In more detail, a standard $\ell_2$ cost is considered for image reconstruction (e.g., deconvolution/super-resolution, semi-blind deconvolution) problems, while a smoothed $\ell_1$ is employed to assess localization precision in some exemplary fluorescence microscopy problems exploiting single-molecule activation. Several numerical experiments are reported to validate the proposed approach on synthetic and realistic ISBI data.
\end{abstract}

\section{Introduction}   \label{sec:intro}

Imaging ill-posed inverse problems are ubiquitous in many applied fields where the quantities observed vary from cells to galaxies  \cite{Bertero-etal-2009}. In the field of fluorescence microscopy, for instance, common sources of degradation are light diffraction and noise interference due to photon-counting processes and the electronics of the acquisition device. Due to these degradation effects, observed data typically reveal little amount of information, being limited in spatial resolution and presenting several artefacts caused by possible optical aberrations and/or background (or auto-)fluorescence effects. Mathematically, the image formation model mapping the unknown informative image $u^{\text{true}}\in\mathbb{R}^n$ onto its corresponding noisy and low-resolution version $f\in\mathbb{R}^m, m<n$ can be naturally formulated as:
\begin{equation}\label{eq:ImageFormation}
	f = N(Au^{\text{true}} + \text{b}),
\end{equation}
where $A=SH\in\mathbb{R}^{m\times n}$ is the product of a convolution matrix $H\in\mathbb{R}^{n\times n}$ describing the convolutional action of the Point Spread Function (PSF) of the instrument on $u$ and $S\in\mathbb{R}^{m\times n}$ is a down-sampling operator modelling the loss of resolution, $b\in\mathbb{R}^m$ is the (possibly space-variant) background image containing out-of-focus fluorescent molecules, while $N:\mathbb{R}^m\to \mathbb{R}^m$ models an interference process introducing noise. In the case of additive white Gaussian noise, for instance, $N(z) = z+n$ with $n\sim \mathcal{N}(0,\sigma^2 \text{Id})$ while in the case of signal-dependent Poisson noise $N(z)=\text{Poiss}(z)$, that is the realisation of a multi-variate Poisson process with mean and variance $z\in\mathbb{R}^m$. Those noise distributions are the ones classically considered in a microscopy setting \cite{ROF,BerteroBook}, although more relevant mixed \cite{Luisier2011,Jezierska2014,CalatroniSIAM2017,Toader2022} and/or multiplicative processes can also be considered. 

\paragraph{Variational regularization: modelling and optimization.} To counteract the instabilities arising when attempting to solve \eqref{eq:ImageFormation} by direct inversion, the field of Bayesian approaches and variational regularisation has proved effective in a plethora of applications, see, e.g., \cite{Stuart_2010,ChanShen} for surveys. In a Bayesian setting, the reconstruction procedure is formulated as an optimisation process in the form
\begin{equation}
\underset{u}{\argmin} \; 
	\left( D(u;f,A,\theta_1) + R(u;\theta_2) =: E_f(u;\theta)\right),  \label{eq:var_form}
\end{equation}
where $D(u;f,A,\theta_1)$ is a data fidelity term whose expression depends on the statistics of the noise, while $R(u;\theta_2)$ is a regularization term which models \emph{a-priori} knowledge on the solutions (such as sparsity or smoothness).
The hyperparameters $\theta \in \R^\Theta$ represents all the (hyper)parameters the energy depends on. The first and most natural of them is the regularization parameter, which we indicate with $\rho \in \R$, and balances the tradeoff between the two terms of the energy.

From an optimization perspective, solving the composite problem \eqref{eq:var_form} may be quite challenging. Convex data and regularization terms (such as, e.g., $D(u;f,A,\theta_1) = \frac{1}{2}\|Au-f\|^2$ and $R(u;\theta_2)=\theta_2\|u\|_1$) have been extensively employed over the last years, motivated by the success of compressed sensing in this field \cite{Candes2006,Donoho2006}. Standard first-order algorithms solving such composite problem (possibly incorporating further convex constraints) rely on the use of proximal-based algorithms. Among them, iterative proximal-gradient \cite{CombettesWajs2005}, primal-dual \cite{ChambollePock2011} and splitting-based approaches \cite{Boyd2011} have been successfully employed (see also \cite{ChambollePock2016} for a review), possibly coupled with acceleration techniques improving the convergence speed \cite{Nesterov-2004,Beck-Teboulle2009}. Whenever proximal points can be expressed in closed-form (as it happens in the case of $\ell_1$ regularization, possibly combined with orthogonal transformations), proximal-gradient or forward-backward algorithms are often regarded as the simplest approaches for solving structured, convex composite problems in the form \eqref{eq:var_form}. For standard $\ell_2$-$\ell_1$ problems the use of such algorithm dates back to \cite{daubechies2004iterative} under the name of Iterative Soft-Thresholding Algorithm (ISTA) which was used as an effective regularization procedures for linear ill-posed inverse problems with sparsity constraints. Given an initialization point $u^0\in\mathbb{R}^n$ and a step-size parameter $\alpha\leq 1/\| A\|^2$ guaranteeing convergence, the ISTA iteration defined for $k\geq 0$ reads:
\begin{equation}  \label{eq:ISTA}
	u^{(k+1)} = \mathcal{T}_{\alpha \theta_2}\left(u^{(k)}-\alpha A^T(Au^{(k)} - f) \right) = \mathcal{T}_{\alpha \theta_2} \left( \left(\text{Id}-\alpha A^T A \right)u^{(k)} + \alpha A^T f\right),
\end{equation}
where $\mathcal{T}_{\alpha\theta_2}:\R^n\to\R^n$ is the soft-thresholding operator defined component-wise by $\mathcal{T}_{\alpha\theta_2}(z)= \max\left( |z|-\alpha\theta_2,0 \right)\text{sign}(z)$.

\paragraph{Parameter estimation.} A crucial challenge in the design of variational models in the form \eqref{eq:var_form} is the choice of the fidelity/regularization parameters $\theta=(\theta_1,\theta_2)$ which, potentially, may lie in a very large dimensional space $\R^{\Theta}$, $|\Theta|\gg 1$. They make the regularization functional $R(\cdot;\theta_2)$ expressive enough to encode the desired solution properties and balance its action against  the data fit which, depending on $\theta_1\in \R^{\Theta_1}$ may encode possible local noise dependence. Classical approaches addressing hyper-parameter estimation problems are based on prior knowledge of the noise level \cite{Morozov1966,Bertero_2010}, (heuristic) study of regularization paths/Pareto fronts and cross-validation, see, e.g. \cite{Hansen1992}. Over the last decade, however, the increasing availability and access to data, favoured the development of machine- and deep-learning based approaches for parameter estimation. The general idea of such methods consist in optimizing model parameters (even when the regularizer is parametrized as a neural-network)  by minimizing some task-dependent quality metric assessing proximity (in some sense) to reference image data, see \cite{Lucas2018,arridge2019,Monga2021} for a review. Having in mind a close connection with the optimization framework of optimization-based iterative schemes such as \eqref{eq:ISTA}, among the plethora of data-driven approaches for imaging, we focus here on the two particular classes of bilevel optimization and algorithmic unrolling techniques, see \cite{SIG-111} and \cite{Monga2021} for review papers, respectively.

\paragraph{Bilevel optimization and algorithmic unrolling.} As it will be described in more detail in Section \ref{sec:bil_unrol}, bilevel learning approaches compute optimal hyper-parameters $\theta$ in \eqref{eq:var_form} by comparing the solution $u^*(\theta)$ of the problem with ground truth data  w.r.t.~some assessment metric. They thus naturally read as a nested optimization problem where the variational model serves to constrain the space of possible solutions. Bilevel approaches have been widely used in the context of variational image reconstruction  \cite{JCCarolaBilevel2013,Kunisch2013,Ochs2016bil,CalatroniBilevel2017,Helsinki,SIG-111} and machine learning \cite{pedregosa16,ghadimi2018approximation,Grazzi2021,NEURIPS2022_aa84ec1a} and offer an interesting and sound perspective for learning hyperparameters in regularized inverse problem formulations.  \emph{Algorithmic unrolling} is a somehow similar approach where an explicit writing of the first $K\in \mathbb{N}$ iterations of an iterative algorithm as a function of the initial point $u^0$ and the model parameters is interpreted as the output $u^{(k)}(\theta)$ of a $K$-layer neural network (with suitable activation functions) to be matched with ground-truth data. Taking the ISTA iteration \eqref{eq:ISTA} as a toy example, in \cite{GregorLeCun2010} an unrolled version of ISTA (therein called Learned ISTA, LISTA) was considered. Using the notation above, the main idea there consists in learning parameters $\left\{\alpha,\theta_2,W_e,W_t\right\}$ such that for $k=0,\ldots, K$
\begin{equation}   \label{eq:LISTA}
    u^{(k+1)}  = \mathcal{T}_{\alpha \theta_2} \left( W_e u^{(k)} + W_t f\right),
\end{equation}
thus resembling \eqref{eq:ISTA} but being potentially more expressive being \emph{a-priori} not connected to the minimization of any underlying variational problem.

Algorithmic unrolling is nowadays a popular strategy to solve inverse problems, see, e.g., \cite{hershey2014deep,Adler2018,Bertocchi_2020,Monga2021} for some relevant references and \cite{Brauer2022,riccio2024} for connections with bilevel optimization and Deep Equilibrium Models. Depending on the particular problem and algorithm at hand, choosing to learn \emph{all} problem parameters at once in an end-to-end fashion may be not ideal especially in applications involving known physical models (typically, encoded in the choice of the operator $A$ or in the modelling of the data-fit $D$) which one would like to keep fixed.

\paragraph{A case study: reconstruction and localization in fluorescence microscopy.}

In fluorescence microscopy imaging, the physical limitations imposed by light diffraction makes the accurate reconstruction and analysis of small biological samples very challenging. For standard microscopes structures closer to $\sim 250$ nm in the $x$-$y$ plane cannot be distinguished. Super-resolution fluorescence microscopy techniques \cite{Dickson1997,Hell94,Betzig2006} aim to overcome such barrier, allowing resolution (i.e., minimal distance between two distinct objects) up to $20$ nm \cite{3D_STORM,iPALM}.
In this context, photons emitted by photo-activable fluorescent molecules are captured by detectors after passing through special lenses and optical devices, which limit spatial resolution due to light diffraction. Such process can be modelled by \eqref{eq:ImageFormation} as a convolution of the emitters with the microscope PSF, thus producing blurred and low-resolution measurements featuring also noise distortions due to electronic interference and a (possibly space-variant) background term accounting for out-of-focus molecules.

Among the many super-resolution techniques proposed in the microscopy community to address this challenge, some have attracted the attention of the applied mathematics community working on signal and image processing. The former has been popularized under the name of Single Molecule Localization Microscopy (SMLM) approaches, see, e.g., \cite{Lelek2021} for a review. Here, the idea consists in acquiring a temporal sequence of images $\left\{ f_t\right\}_{t=1}^T$ where at each frame $t$ only a small percentage of special photo-activable fluorescent molecules is active, thus making the detection of nonzero elements in the corresponding image $u_t$ (that is, the localization process) easier. A super-resolved image $u_{\text{SR}}$ is thus obtained by averaging $u_{\text{SR}}=\frac{1}{T}\sum_{t=1}^T u_t$. A mathematical modelling of this problem relies, essentially, on the use of sparsity-based regularization techniques for accurate localization, either in a convex \cite{Aritake2021} or nonconvex \cite{Gazagnes2017,LazzarettiISBI2021} regime. Parameter selection plays here a crucial role and their brute-force optimization could be a tedious task.  Other approaches to SMLM make use of end-to-end deep learning procedures for accurate localization upon suitable training \cite{Nehme_18}. 

Some other techniques require less specific fluorescent molecules being based on the estimation of a super-resolved image by exploiting only a temporal sequence of images whose intensity fluctuations are analyzed typically in terms of their second-order statistics \cite{SOFI,SPARCOM,srrf,Stergiopoulou2021,stergiopoulou2022}. Those techniques are better suited to standard biological set ups where the repeated ON-OFF processes typical of SMLM may be harmful for the sample under observation. In these settings, a modelling similar to \eqref{eq:ImageFormation} but reformulated in a covariance domain is considered (see Section \ref{subsec:Exp4} for more details).

\paragraph{Contribution.} In this work, we unroll an accelerated projected gradient descent (APGD) scheme to estimate optimal parameters for a (smoothed) $\ell_2$-$\ell_1$ regularization problem endowed with a non-negativity constraint. In the framework of fluorescence microscopy imaging, this simple regularization model can be effectively used for both reconstruction and molecule localization purposes. Given the intrinsic difference between these two tasks, we consider in the following task-adapted evaluation metrics estimating optimal algorithmic parameters depending on the specific objective considered. We report several numerical experiments confirming the validity of the approach on exemplary deconvolution/super-resolution problems classically encountered in fluorescence microscopy, also in a semi-blind scenario in the context of fluctuation-based microscopy for estimating the shape of the blurring function encoded by the operator $A$.

\section{Models and algorithms}  \label{sec:bil_unrol}

We review the mathematical formulation of the bilevel optimization schemes considered and consider their unrolled versions, making precise the models considered to deal with the reconstruction/localization problems encountered in the microscopy applications considered.

\subsection{From bilevel optimization to algorithmic unrolling}
Bilevel optimization approaches \cite{JCCarolaBilevel2013,Kunisch2013,Ochs2016bil,pedregosa16,CalatroniBilevel2017,ghadimi2018approximation,Grazzi2021,NEURIPS2022_aa84ec1a,Helsinki,SIG-111}  rely on a supervised data-driven approach to estimate model and algorithmic hyperparameters. Given a training dataset $\cD = \{(f_t,g_t) \colon t = 1,\ldots,T\}$, where $f_t$ is a corrupted (noisy, blurred, low-resolution) version of the ground truth image  $g_t$, the problem problem of estimating optimal parameters $\widehat{\theta}$ can be formulated
\begin{equation}\label{eq:BilevelProblem}
	\begin{cases}
		\widehat{\theta}\in \underset{\theta \in \R^p}{\argmin} \; \sum_{t = 1}^T \cL(u^*_t(\theta),g_t) \\
		\mathrm{s.t.} \quad  u^*_t(\theta) \in \underset{u}{\argmin} \; E_{f_t}(u;\theta) \qquad t = 1,\ldots,T.                 
	\end{cases}
\end{equation}
Essentially, it consists in two nested optimization problems, where the inner (or lower) one feeds its result to the outer (or upper) one which is used to assess optimality of the estimation w.r.t.~task-dependent evaluation loss. The ultimate goal of the scheme is to find an optimal set of parameters $\widehat{\theta}$ which optimizes the performance of the underlying variational model. 

Solving \eqref{eq:BilevelProblem} requires the development of (gradient-based) optimization methods searching for the optimal parameters configuration minimizing the outer loss function. For that, the computation of the minimizers $u_t^*(\theta)$ for $t=1,\ldots,T$ of lower-level functional is required,  which is often unfeasible. Most of the time, it is thus necessary to use iterative schemes solving the variational lower-level constraint under a fixed computational budget (that is fixed amount of iterations), thus leading to approximate minimizers and, ultimately,  inexact computations of the gradient of the outer loss. This idea stands at the very basis of algorithmic unfolding \cite{hershey2014deep,Monga2021}. Mathematically this can be modelled by replacing \eqref{eq:BilevelProblem} with the unrolling over $K\geq 1$ iterations of an optimization algorithm  $\mathcal{A}$, such that:
\begin{equation}\label{eq:Unrolling}
    u_t^*(\theta) \approx \cA^K(f_t;\theta)
\end{equation}
so that optimization is performed for a fixed number of iterations $K$. By expanding the expression of the computed quantity from $k=K$ till $k=0$ it is thus possible to have a compact expression of $u^K_t(\theta)$ which can be used for differentiating over the desired hyperparameters $\theta$, see, e.g., \cite{Ochs2016bil}. Note that, depending on the particular choice of $\mathcal{A}$, in addition to model parameters (such as regularization parameters and/or quantities related to the forward model), algorithmic parameters (such as algorithmic step-sizes) can also be learned. Once the training phase has been completed and an optimal configuration of parameters $\widehat{\theta}$ has been found, one can simply apply the trained algorithm $\cA^K(\cdot;\widehat{\theta})$ on new unseen data.

\subsection{A model for sparse reconstruction/localization} \label{sec:modelSparse}

We describe in this section the variational model employed to solve localization and reconstruction problems often encountered in the framework of microscopy as in the case, for instance, of Single Molecule Localization methods \cite{Aritake2021,Gazagnes2017,LazzarettiISBI2021}. In order to consider an unrolled bilevel strategy \eqref{eq:BilevelProblem}--\eqref{eq:Unrolling}, some details on the the algorithm $\mathcal{A}$ considered and on the training procedure employed (with a particular focus on the choice of the loss function $\mathcal{L}$ are also given. 

\paragraph{Localization and reconstruction by sparse optimization.}
We consider the following $\ell_1$ sparse reconstruction approach with a non-negativity constraint:
\begin{equation}\label{eq:E}
E_f(u;\theta) = D(Au;f,bg) + \rho \|u\|_1 + \iota_{\R^n_+}(u), \qquad \iota_{\R^n_+}(u) := \begin{cases}
    0       & u \in \R^n_+ \\
    +\infty & \text{otherwise}
    \end{cases},
\end{equation}
where $\R^n_+ = \{u \in \R^n \colon u_i \geq 0, i = 1,\ldots,n\}$ and $bg \in \R$ represents a constant positive background term which can be used both for modelling, see, e.g., \cite{Harmany2012}. Note that the non-negativity constraint makes the $\ell_1$-norm differentiable. In this work we make use of the $\ell_1$-norm as a convex approximation of the $\ell_0$ pseudo-norm in a compressed sensing fashion \cite{Candes2006,Donoho2006}. 

The fidelity term $D$ enforces proximity between the physical modelling and the observed data depending on the particular modelling of the noise statistics.  In particular we consider the two choices:

\begin{enumerate}[(i)]
    \item (Additive white Gaussian noise) $D(Au;f) = \frac{1}{2} \|Au - f\|_2^2$
    \item (Poisson noise) $D(Au;f,bg) = KL(Au+bg;f) = \sum_{i = 1}^m (Au)_i + bg - f_i - f_i \log\left( \frac{(Au)_i + bg}{f_i}\right)$.
\end{enumerate}

Fidelity (i) is the quadratic penalty classically employed to model the presence of additive white Gaussian noise \cite{ROF}, while fidelity (ii) is the Kullback-Leibler (KL) used to model signal-dependent Poisson noise \cite{Bertero2008,Harmany2012}. Note that both functionals are convex and Lipschitz differentiable on the non-negative orthant. However, the KL divergence might present some numerical instability issues due to the presence of the logarithm for very small values of original ($f_i$) and/or reconstructed $(Au)_i + bg$ components of the data. 
We consider in the following an algorithmic unfolding procedure estimating optimal hyperparameters for \eqref{eq:E} coupled for the two data terms above to assess and compare the performance of both models on localization and reconstruction problems.

\paragraph{Solving the lower-level problem via accelerated projected gradient descent.}
By solving \eqref{eq:BilevelProblem} via algorithmic unrolling \eqref{eq:Unrolling}, the computation of the loss gradient can simply be done via backpropagation, provided that the single iteration of $\cA$ is differentiable. Given the composite structure of problem \eqref{eq:E}, we consider as $\cA^K(f;\theta)$ $K$ iterations of an accelerated projected gradient descent endowed with a smoothed differentiable projection initialized by taking $f$ as starting point:
\begin{algorithm}
	\caption{Accelerated Projected Gradient Descent}\label{alg:FISTA}
	\textsc{Input}: $f \in \R^m$.\\
	$
	\begin{array}{ll}
		u^{(0)}   = u^{(-1)}  = f\\
		\mbox{\textsc{For} }  k = 0,...,K-1\\
		\left\lfloor\begin{array}{lcl}
		\bar v^{(k)} & = & u^{(k)} + \beta_k(u^{(k)}-u^{(k-1)})          \\[0.2cm]
		v^{(k)}      & = & \Pi(\bar v^{(k)})                             \\[0.2cm]
		w^{(k)}      & = & v^{(k)}-\alpha_k \nabla_u E_{f}(v^{(k)};\theta) \\[0.2cm]
		u^{(k+1)}    & = & \Pi(w^{(k)})                                  
	\end{array}\right.\\
	u^* = u^{(K)}
	\end{array}
	$\\
	\textsc{Output}: $u^*$.
\end{algorithm}

 The map $\Pi: \R^n \rightarrow \R_+^n$ has been used in  \cite{Helsinki}: it is the smoothed version of the standard (non-differentiable) Euclidean projection. In particular, for an arbitrarily small $\varepsilon>0$, $\Pi$ acts component-wise on $u\in\R^n$ as:
\[
\left(\Pi(u)\right)_i = \begin{cases}
                        \left(2 - \frac{u_i}{\varepsilon}\right)\frac{u_i^2}{\varepsilon},\quad &
                    \quad u_i \in (0,\varepsilon) \\
                        \text{max}\left(u_i,0\right)\quad &\quad \text{otherwise}
                        \end{cases}.
\]
By definition $\Pi$ acts as the Euclidean projection outside the interval $(0,\varepsilon)$, while it coincides with a third degree polynomial interpolating the points $(0,0)$ and $(\varepsilon,\varepsilon)$.

Algorithm \ref{alg:FISTA} is an instance of Accelerated Projected Gradient descent \cite{Nesterov-2004,Beck-Teboulle2009} with a smoothed projection. The sequence $\{\beta_k\}_{k=0}^{K-1}$ of inertial parameters can be  chosen as 
\[
\beta_k = \begin{cases}
            0  & \text{if } k = 0 \\
            \frac{k-1}{k+2} & \text{if } k \geq 1
          \end{cases},
\]
which was shown in  \cite{Chambolle_2015_FISTA_Iterates} to guarantee weak convergence of the iterates. The sequence of step-size parameters $\{\alpha_k\}_{k=0}^{K-1}$ can be set either to a constant equal to the inverse of the Lipschitz constant of $\nabla_u E_f$, or adaptively chosen using, e.g., monotone Armijo-type rules \cite{Beck-Teboulle2009}. In the following we let this sequence to be learned.

\medskip

 Since the general $k$-th iteration of $\cA$ is differentiable, it is possible to compute the gradient of a differentiable loss function by simple backpropagation procedure. Namely, for one training sample $(f,g)$ and one parameter $\theta_j$ we have:
\begin{align}
    \frac{\partial \cL(u^{(K)}(\theta),g)}{\partial \theta_j} & = \left( \frac{\partial u^{(K)}(\theta)}{\partial \theta_j}\right)^T\frac{\partial\cL(u^{(K)}(\theta),g)}{\partial u}.
\end{align}
The chain rule can then be applied through $\frac{\partial u^{(K)}(\theta)}{\partial \theta_j} = \frac{\partial}{\partial \theta_j}\cA(u^{(K-1)};\theta)$. With a small abuse of notation and denoting by $\Pi'(v)$ the diagonal Jacobian matrix of $\Pi$ evaluated at $v$, we have
\begin{align}
    \frac{\partial u^{(K)}(\theta)}{\partial \theta_j} & = \Pi'(w^{(K-1)})\frac{\partial}{\partial \theta_j}\left(v^{(K-1)} - \alpha_{K-1}\nabla_uE_f(v^{(K-1)};\theta)\right) \\
     & = \Pi'(w^{(K-1)})\left(I_n - \alpha_{K-1}\nabla^2_{uu}E_f(v^{(K-1)};\theta) \right)\frac{\partial v^{(K-1)}}{\partial \theta_j} \\                               & \hspace{12pt} - \alpha_{K-1}\Pi'(w^{(K-1)})\frac{\partial}{\partial \theta_j}\nabla_u E_f( v^{(K-1)};\theta),
\end{align}
where $I_n$ is the identity matrix of size $n\times n$ and 
\begin{equation}
    \frac{\partial v^{(K-1)}}{\partial \theta_j} = \Pi'(\bar{v}^{(K-1)})\left((1+\beta_{K-1})\frac{\partial u^{(K-1)}(\theta)}{\partial \theta_j} - \beta_{K-1} \frac{\partial u^{(K-2)}(\theta)}{\partial \theta_j}\right),
\end{equation}
so that the computations can be repeated for $\frac{\partial u^{(K-1)}(\theta)}{\partial \theta_j}$ and $\frac{\partial u^{(K-2)}(\theta)}{\partial \theta_j}$ until $u^{(0)} = f$.

In the computations above, the quantities $\nabla_{uu}^2E_f(v^{(k)};\theta)$ and $\nabla^2_{\theta u}E_f(v^{(k)};\theta)$ are present, which requires the conditions:
\begin{itemize}
    \item[-] $E$ is twice continuously differentiable, with $\nabla_{uu}^2E$ invertible;
    \item[-] $\nabla_u E$ is differentiable with respect to $\theta$.
\end{itemize}
Both properties are satisfied by the energy function \eqref{eq:E}. Note that such conditions are quite restrictive for more effective sparsity-promoting regularization terms approximating the $\ell_0$ pseudo-norm, such the continuous exact relaxations CEL0 studied in \cite{Soubies_2015_CEL0,Soubies2017}. Indeed, just like the $\ell_1$-norm, such penalty is differentiable with respect to $u$ when the non-negativity constraint is taken into account, but its gradient is not differentiable with respect to the regularization parameter that governs it (see Appendix \ref{app:A}). 

\paragraph{The training loss functions.}
In order to adapt the hyperparameter optimization to the particular task at hand (reconstruction/localization), we consider in the following experiments two different assessment metrics: 
\begin{itemize}
    \item (Reconstruction loss)  a standard $\ell_2$-norm (corresponding to SNR optimization, as done, e.g., in \cite{JCCarolaBilevel2013,Kunisch2013}):
    \begin{equation} 
       \mathcal{L}_2(u^*(\theta);g_t) = \frac{1}{2} \|u^*(\theta) - g_t\|^2_2, \label{eq:LossL2}
       \end{equation}
    \item (Localization loss) a smoothed $\ell_1$-norm with image binarization assessing localization precision:
    \begin{equation}
        \mathcal{L}_1(u^*(\theta);g_t)  = \sum_{i=1}^{n}\psi_{\gamma}\left(\left(B_{\delta,c,\epsilon}(u^*(\theta)) - \tilde{g_t}\right)_i^2\right). \label{eq:LossL1} 
    \end{equation}
\end{itemize}

Note that in (\ref{eq:LossL2}) $g_t$ is a ground-truth image from the training dataset while in (\ref{eq:LossL1}) $\tilde{g_t}$ is a binarized version of $g_t$ where each nonzero element is set as 1 or 255, depending on the scaling. As observed in \cite{Sage2019}, in the field of microscopy (and in particular of single-molecule localization approaches) is important to define tailored assessment metrics adapted to the task at hand which may aim at either estimating the precise intensity of the light they emit (and for such case \eqref{eq:LossL2} is a natural choice) or detect whether model \eqref{eq:E} has computed a reconstruction in the correct position, for which a loss in the form \eqref{eq:LossL1} seems more natural. To allow the computation of derivatives, we employ in  definition   \eqref{eq:LossL1} a function $\psi_{\gamma}$ which is the Huber function defined, for $\gamma>0$, by:
\begin{equation}\label{eq:Huber}
    \psi_{\gamma}(s) := \begin{cases}
                            \frac{s}{\gamma}  & s \leq \gamma^2 \\
                            2\sqrt{s} - \gamma & s^2 > \gamma^2.
                        \end{cases}
\end{equation}

\begin{figure}[!h]
    \centering
    \includegraphics[scale = 0.35]{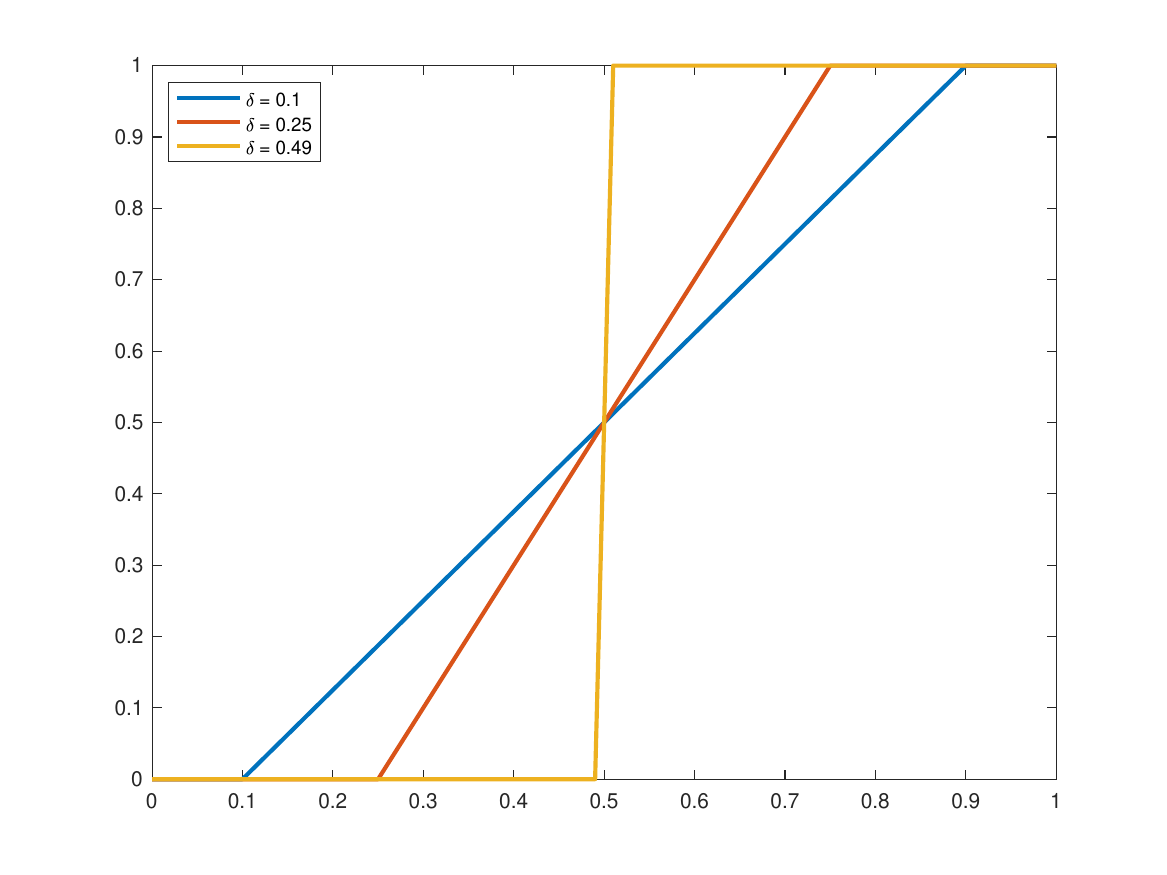}
    \caption{Binarization in [0,1] with $c = 0.5$, $\epsilon = 10^{-4}$ for three different values of $\delta$: 0.1 (blue), 0.25 (red), 0.49 (yellow).}
    \label{fig:bin}
\end{figure}

The function $B_{\delta,c,\epsilon}:\R^n\to\R^n$ acts pixel-wise by binarizing the tentative reconstruction $u^*_t(\theta)$ for the current estimation of $\theta$ in order to compare it with $\tilde{g}_t$. For a given image $s\in\R^n$, by setting $\bar{p} 
 =\max_i~s_i$, it is defined as
\begin{equation}\label{eq:Bin}
    B_{\delta,c,\varepsilon}(s) := 
        \begin{cases}
            0 & s \leq \delta, \\
            \left( 2 - \frac{s - \delta}{\varepsilon} \right) 
            \frac{(s-\delta)^2}{\varepsilon} \frac{\bar p}{2(c-\delta)} & \delta < s < \delta + \varepsilon, \\
            \frac{\bar p}{2(c-\delta)} (s - \delta) & \delta + \varepsilon \leq s \leq 2c - \delta - \varepsilon, \\
            \bar p - \left( 2 - \frac{2(c - \delta) - s + \delta}{\varepsilon} \right) \frac{(2(c - \delta) - s + \delta)^2}{\varepsilon}\frac{\bar p}{2(c-\delta)} & 2c - \delta - \varepsilon < s < 2c - \delta, \\
            \bar p & s \geq 2c - \delta.
        \end{cases} 
\end{equation}
We observe that such definition depends on three parameters:  $\delta>0$ is the thresholding parameter used to binarize into zeros and the value $\bar p$, so that the point $s = c$ is the one getting transformed into $\frac{\bar p}{2}$. The parameter $\varepsilon$ controls the smoothing intervals, which are symmetric w.r.t.~$c$, in the same way as in the projection $\Pi$ previously described. Figure \ref{fig:bin} shows various shapes of $B_{\delta,c,\epsilon}$ for three different values of $\delta$. Ideally, we would like to consider situations where $c\approx \delta$ to have a sharp transition between the two different regimes. 
Note that while this choice may look at a first sigh a complicated way of assessing localization precision (as it introduces further parameters to estimate), in practice  $\epsilon$ does not need any fine-tuning provided it is sufficiently small. To obtain a steep slope one could further choose $c = \delta + c_0$, with $c_0 > \epsilon$, still small. Under this choice, the binarization actually becomes differentiable with respect to $\delta$, which allows for a strategy to learn $\delta$ as well. There holds indeed:
\begin{equation}
\frac{\partial \mathcal{L}_1}{\partial \delta} = \sum_{i = 1}^n \psi_{\gamma}'\left((B_{\delta,c,\epsilon}(u^*(\theta)) - \tilde{g_t})_i\right) \frac{\partial B_{\delta,c,\epsilon}(u^*(\theta))_i}{\partial \delta}.
\end{equation}
where
\begin{equation}
\frac{\partial B_{\delta,c,\epsilon}(s)}{\partial \delta} = 
\begin{cases}
    0                                                                & \delta \geq s \\
    \frac{\bar p}{2c_0} \left(3 \frac{(s-\delta)^2}{\epsilon^2} - 4 \frac{(s-\delta)}{\epsilon}\right)                                               & s > \delta > s - \epsilon  \\
    -\frac{\bar p}{2 c_0}                                            & s - \epsilon \geq \delta \geq s + \epsilon - 2 c_0 \\
    -\frac{\bar p}{2c_0}  \left(4 \frac{2c_0 - s + \delta}{\epsilon} - 3 \frac{(2c_0 - s + \delta)^2}{\epsilon^2}\right)                       & s + \epsilon - 2 c_0 > \delta > s - c_0 \\
    0                                                                & \delta \leq s - 2c_0
\end{cases},
\end{equation}
which can be used for defining gradient-type updates.

\section{Numerical experiments} 

In this Section we illustrate the results obtained by unrolling Algorithm \ref{alg:FISTA} applied to model \eqref{eq:E} for exemplar reconstruction/localization  problems using different modelling for the noise statistics and loss functions \eqref{eq:LossL2}--\eqref{eq:LossL1} for training the model.

\begin{figure}
    \centering
    \begin{subfigure}{.475\linewidth}
    \centering
        \includegraphics[scale = 0.15]{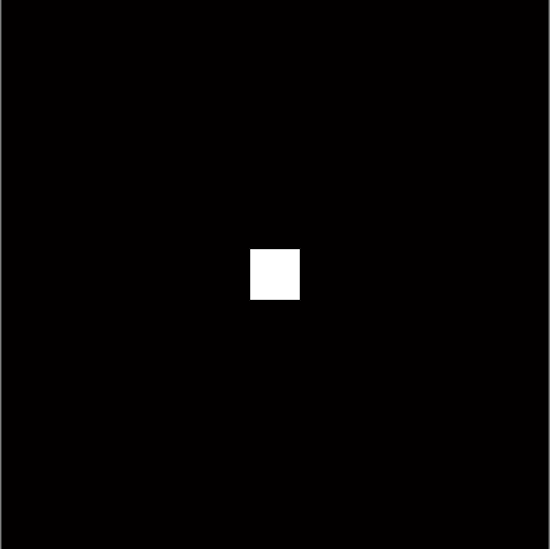}
        \caption{}
        \label{fig:jaccard_example_subfig1}
    \end{subfigure}
    \begin{subfigure}{.475\linewidth}
        \centering
        \includegraphics[scale = 0.15]{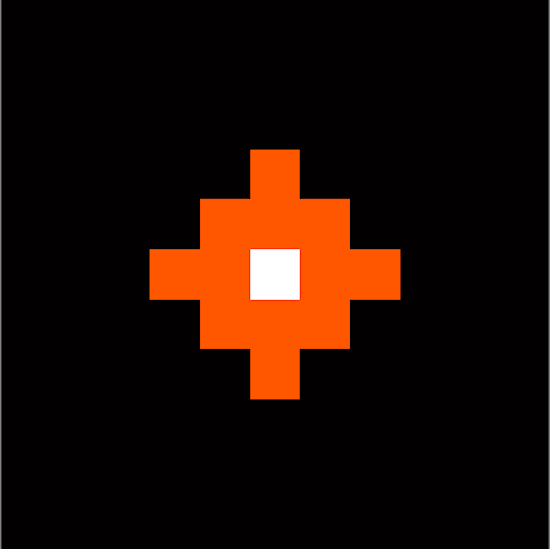}
        \caption{}
        \label{fig:jaccard_example_subfig2}
    \end{subfigure}
    \\
    \begin{subfigure}{.475\linewidth}
        \centering
        \includegraphics[scale = 0.15]{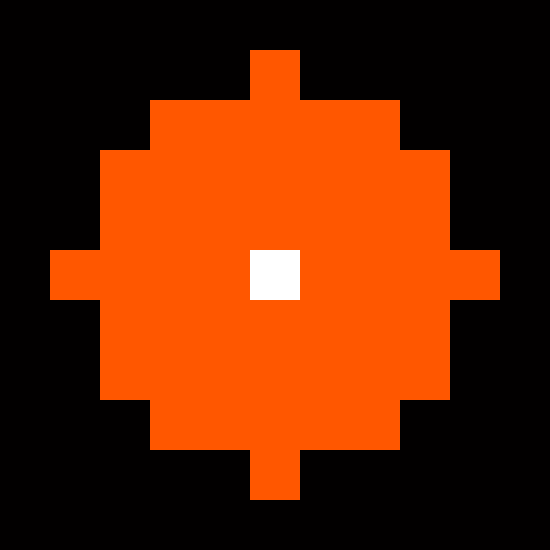}
        \caption{}
        \label{fig:jaccard_example_subfig3}
    \end{subfigure}
    \begin{subfigure}{.475\linewidth}
        \centering
        \includegraphics[scale = 0.15]{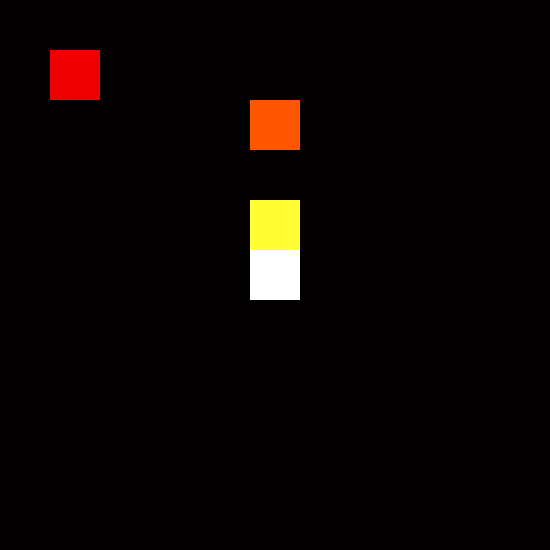}
        \caption{}
        \label{fig:jaccard_example_subfig4}
    \end{subfigure}
    \caption{(a) Radius of correct detection for $\tilde{\delta} = 0$. (b) Radius of correct detection for $\tilde{\delta} = 2$. (c) Radius of correct detection for $\tilde{\delta} = 4$. (d) Example in the case of $\tilde{\delta} = 0$: the white pixel is the true molecule, the yellow pixel is counted as a TP, the red pixel is FP as it is outside of the range, while the orange pixel is still counted as a FP even if it is in range because the yellow one is closer to the true molecule.}
    \label{fig:jaccard_example}
\end{figure}

Unless specified otherwise, we evaluate the quality of the images by averaging, over all the $T$ frame frames, the Peak Signal to Noise Ratio (PSNR) and the Jaccard index defined by:
\begin{equation}\label{eq:Jaccard}
J_{\tilde{\delta}} = \frac{\# \text{ TP}}{(\# \text{ TP}) + (\# \text{FN}) + (\# \text{FP})} \in [0,1]
\end{equation}
where TP, FN and FP denote True Positives, False Negatives and False Positives, respectively. The parameter $\tilde{\delta}>0$ is a tolerance parameter acting as follows:  a reconstructed pixel is counted as TP if it lies within a ball of radius $\tilde{\delta}$ from a true molecule. Figure \ref{fig:jaccard_example} contains a visual representation of how the size of the tolerance value $\tilde{\delta}$ influences the value of the Jaccard index. Given the a ground truth pixel $a\neq 0$, let $\mathcal{B}_{\tilde{\delta}}(a)$ the ball of radius $\tilde{\delta}$ around $a$. A point reconstructed within $\mathcal{B}_{\tilde{\delta}}(a)$ is then counted as TP, outside it is counted as a FP. A TN pixel corresponds to a GT point does not corresponding to any reconstruction within $\mathcal{B}_{\tilde{\delta}}(a)$.  In Figure \ref{fig:jaccard_example_subfig1} ($\tilde{\delta}=0$) only $a$ is counted as a TP, whereas for $\tilde{\delta}=2,4$ (Figures \ref{fig:jaccard_example_subfig2}--\ref{fig:jaccard_example_subfig3}) all points lying within the orange region are counted as TP. Note, however, that if  two elements are computed within $\mathcal{B}_{\tilde{\delta}}(a)$, only the closest to the GT point will be counted as a TP, with the other being a FP, see Figure \ref{fig:jaccard_example_subfig4}. 

In all experiments, the optimization of the upper level problem w.r.t. $\theta$ is carried out via the Scaled Gradient Projection method (SGP) \cite{Bonettini_2009_SGP},
so that we can impose reasonable bounds on the parameters. Some of such bounds are natural (such as the non-negativity of step-size parameters), but for some others some comments are needed. In particular, as far as the choice of the regularization parameters is concerned, in the case of a least squares term, we followed \cite{Koulouri_2021_AdaptiveSuperRes} and chose an upper bound on $\rho$ given by
\begin{equation}
    \rho_{\text{max}} = \min_{t} \|A^Tf_t\|_{\infty}.
\end{equation}
It is easy to check that for all $\rho>\rho_{\text{max}}$ the trivial zero solution is obtained.  For the Poisson data term, such upper bound was chosen empirically.
The upper bound on $\delta$ was set at half the maximum pixel value to avoid null reconstructions. The lower bounds were all set to $10^{-10}$.

Due to the nonconvex nature of the bilevel optimization problem  \eqref{eq:BilevelProblem}, the initial set of parameters $\theta^{(0)}$ may have a major impact on the final results. 
For all set of experiments the initial configurations are reported in the tables with the results. Note that the steplengths $\{\alpha_k\}_{K=0}^{K-1}$ are all initialized to one same value $\alpha$.

\subsection{Single molecule localization microscopy.}\label{subsec:Exp1}

\subsubsection{Simulated data}

We now consider a Single Molecule Localization Microscopy framework where taking as a reference model \eqref{eq:E}, the unrolling strategy is used to learn optimal parameters   $\theta = (\rho,\alpha_0,\ldots,\alpha_{K-1})$ when considering the loss in \eqref{eq:LossL2} and $\theta = (\rho,\delta,\alpha_0, \ldots,\alpha_{K-1})$ for the loss in \eqref{eq:LossL1}.

We first test our strategy on simulated data in order to check whether the Localization loss \eqref{eq:LossL1} promotes better localization properties than \eqref{eq:LossL2}.
We consider a data set of 10 ground truth images $g_t$ containing each between 75 and 150 randomly activated pixels with values in the range $U=[100,255]$. The size of these (vectorized) images is $n\times n = mL \times mL$, where $L = 4$ is the super-resolution factor and $m = 64$ is the size of the downsampled image. The underlying linear model is described by an operator $A = SH$ being the composition of a convolution operator $H \in \R^{n^2 \times n^2}$ describing the convolution of a Gaussian PSF with standard deviation 2.5, and $S \in \R^{m^2 \times n^2}$, a down-sampling operator summing up pixel intensities for each $4\times 4$ patch of the fine grid. Denoting by $\mathbf{1}_L$ and $\mathbf{0}_L$ the $L \times 1$ arrays with constant entries equal to 1 and 0, respectively, and considering the matrix
\begin{equation*}
    S_L =
    \begin{pmatrix}
        \mathbf{1}_L & \mathbf{0}_L & \cdots & \mathbf{0}_L \\
        \mathbf{0}_L & \mathbf{1}_L & \cdots & \mathbf{0}_L \\
        \vdots       &              &        &              \\
        \mathbf{0}_L & \mathbf{0}_L & \cdots & \mathbf{1}_L
    \end{pmatrix} \in \R^{m \times n},
\end{equation*}
then the action of the operator $S$ on a vectorized image  $x$ can be written as
\begin{equation*}
    Sx = \text{vec}\left(S_L\;\bar{x}\;S_L^T\right),
\end{equation*}
where $\text{vec}(\bar{x})=x$ 
The resulting blurred and downsampled images are then corrupted with additive white Gaussian noise with standard deviation $\sigma=0.15$. As for the unrolling, $K = 190$ iterates were used and the value $c=\delta + 0.01$ is used in \eqref{eq:Bin}.

\smallskip

\begin{figure}[!ht]
  \centering
  \begin{subfigure}{.24\linewidth}
    \centering
    \includegraphics[width = 0.85\linewidth]{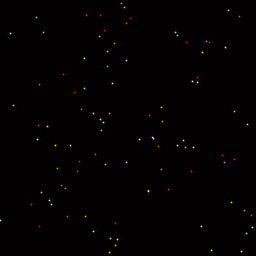}
    \caption{}
  \end{subfigure}%
  \begin{subfigure}{.24\linewidth}
    \centering
    \includegraphics[width = 0.85\linewidth]{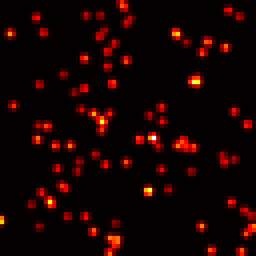}
    \caption{}
  \end{subfigure}%
  \begin{subfigure}{.24\linewidth}
    \centering
    \includegraphics[width = 0.85\linewidth]{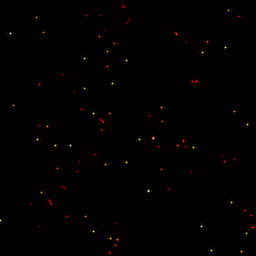}
    \caption{}
  \end{subfigure}
    \begin{subfigure}{.24\linewidth}
    \centering
    \includegraphics[width = 0.85\linewidth]{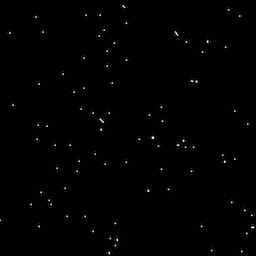}
    \caption{}
  \end{subfigure}
  \caption{(a) Ground Truth $g_t$. (b) Acquired $f_t$. (c) Optimal reconstruction $u(\hat{\theta})$ obtained by optimizing $\mathcal{L}_2$ \eqref{eq:LossL2}. (d) Optimal reconstruction $B_{\hat{\delta},c,\varepsilon}(u^{(K)}(\hat{\theta}))$ obtained by optimizing $\mathcal{L}_1$ \eqref{eq:LossL1}.}
    \label{fig:exp1_recs}
\end{figure}

Table \ref{tab:exp1} shows the values of the Jaccard index \eqref{eq:Jaccard} obtained on the 25 test images for different tolerance values $\tilde{\delta}\in\left\{0,2,4 \right\}$. At a first glance we notice that training the model with \eqref{eq:LossL1} corresponds indeed to an improvement of the Jaccard values w.r.t. training with \eqref{eq:LossL2}. We remark that, in order to obtain meaningful results with \eqref{eq:LossL2}, a tailored dataset had to be considered. Our first tests showed that using, during the training, images with the same density of molecules as those in the test set, the quadratic loss would struggle to find a value of $\rho$ big enough so that the images would not be sufficiently sparse. To avoid that, we created a special dataset of very sparse images containing only 4 non-zero pixels with random values in $U$ and used it to train the model with the $\mathcal{L}_2$ loss. This was done to assess whether reducing the density of the activated pixels would lead to an improvement in terms of Jaccard index. Some exemplary results are reported in Figure \ref{fig:exp1_recs}.

\begin{table}[!h]
    \centering
    \caption{Results for Experiment 1. R = Number of examples used in the training; RS = Number of special examples used in the training. The test set is made by 25 images. For each quality measure, the best score across the models is in bold. The last column reports the initial parameter configuration} 
    \label{tab:exp1}
    \vskip 3mm
    \begin{tabular}{c c | c | c | c | c | c | l}
        \toprule
        Loss function  & R / RS & $J_0$ & $J_2$ & $J_4$ & Avg. J & PSNR & Initial $\rho$, $\alpha$, $\delta$ \\
        \midrule
        $\mathcal{L}_2$  & 0 / 8  & 0.3878 & 0.4179 & 0.4179 & 0.4079 & \textbf{35.76} & 0.1, 0.1, n.a. \\
        \midrule
        $\mathcal{L}_1$  & 10 / 0 & \textbf{0.7452} & \textbf{0.8480} & \textbf{0.8583} & \textbf{0.8172} & 33.09 & 0.1, 0.1, 50 \\
        \bottomrule
    \end{tabular}
\end{table}

\subsubsection{ISBI super-resolution datasets.}\label{subsec:Exp2}

We now test the unrolled localization/reconstruction procedures on the ISBI 2013 SMLM challenge data\footnote{\url{https://srm.epfl.ch/srm/index.php}}. The dataset features 361 ground truth images. Each one has approximately 200 pixels activated, with values in $\{1,2,3\}$. For our experiments we rescaled them so that the maximum value was $255$. The PSF was obtained using information provided by the challenge organizers: the observed data have size $m \times m = 64 \times 64$ pixels of side 100\textit{nm}. The PSF considered is a Gaussian convolution kernel with Full Width at Half Maximum (FWHM) of 258.21\textit{nm}. We consider $L=4$ as a super-resolution factor.

We created a training dataset $\left\{g_t,f_t\right\}$ simulating different levels and types of noise:
\begin{itemize}
    \item[-] Gaussian noise with standard deviation $0.15$.
    \item[-] Poisson noise with constant low background emission $b=0.1$. 
    \item[-] Poisson noise with constant high background emission $b=12.75$.
\end{itemize}

Given the signal-dependence of Poisson noise on the signal, the higher is the background $b$ and the higher is the noise. The two scenarios considered thus model low- and high-photon counts and serve us to assess the performance of the approach over different noise statistics. Reconstruction and localization performances were thus tested in this experiment for both upper-level loss functions \eqref{eq:LossL2}--\eqref{eq:LossL1} and for the Gaussian/Poisson data fidelity functions.

\begin{figure}[!ht]
    \centering
    \begin{subfigure}{0.3\linewidth}
        \centering
        \includegraphics[width = 0.7\textwidth]{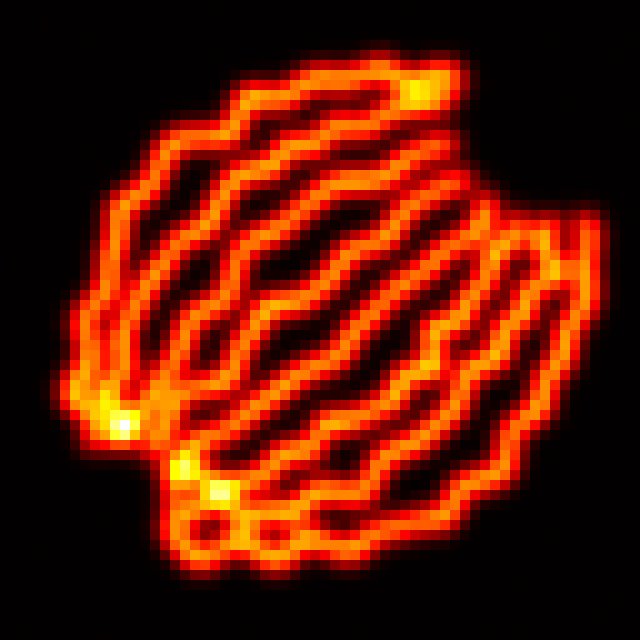} 
        \caption*{}
    \end{subfigure}%
    \begin{subfigure}{0.3\linewidth}
        \centering
        \includegraphics[width = 0.7\textwidth]{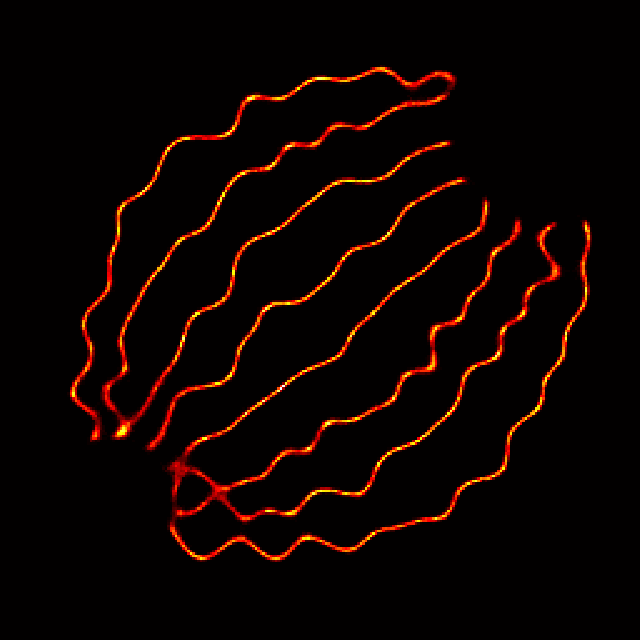} 
        \caption*{$J_4 = 0.7097$}
    \end{subfigure}%
    \begin{subfigure}{0.3\linewidth}
        \centering
        \includegraphics[width = 0.7\textwidth]{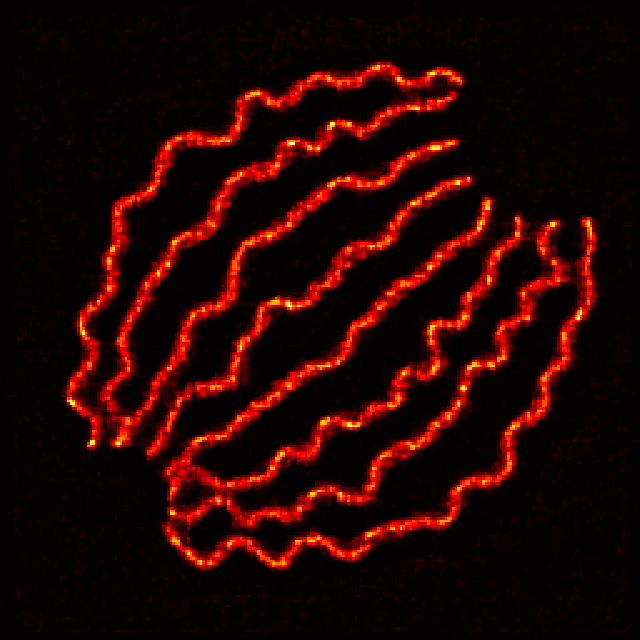} 
        \caption*{$J_4=0.2959$}
    \end{subfigure}\\
        \begin{subfigure}{0.3\linewidth}
        \centering
        \includegraphics[width = 0.7\textwidth]{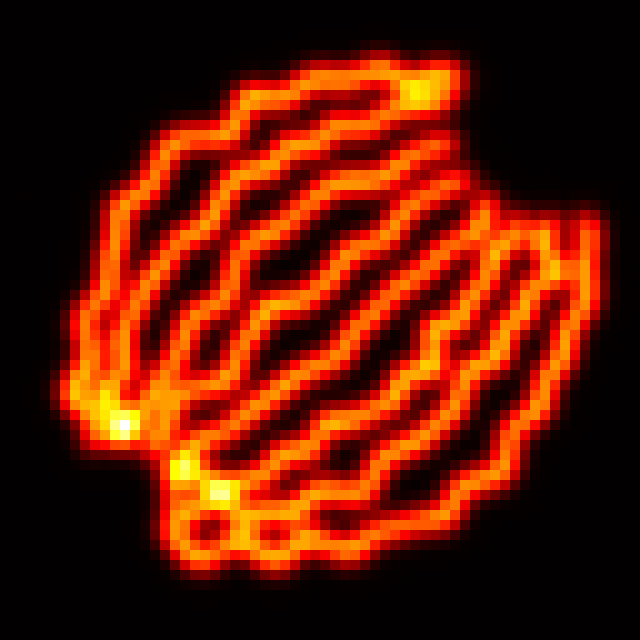} 
        \caption*{}
    \end{subfigure}%
    \begin{subfigure}{0.3\linewidth}
        \centering
        \includegraphics[width = 0.7\textwidth]{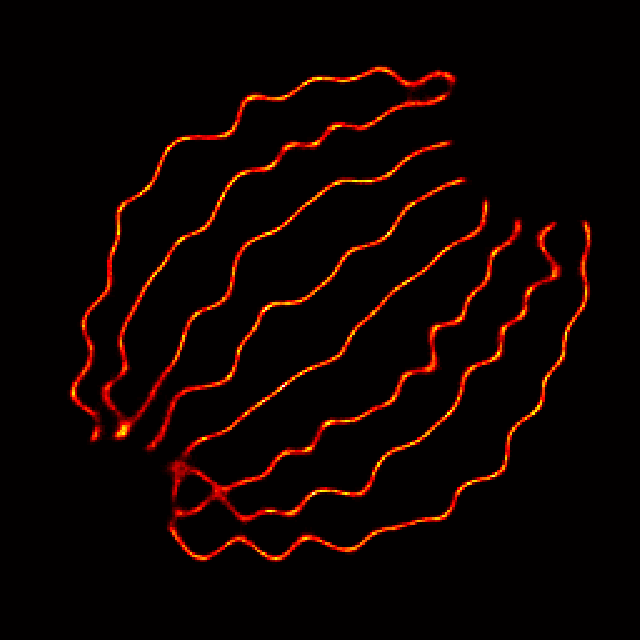} 
        \caption*{$J_4 = 0.6551$}
    \end{subfigure}%
    \begin{subfigure}{0.3\linewidth}
        \centering
        \includegraphics[width = 0.7\textwidth]{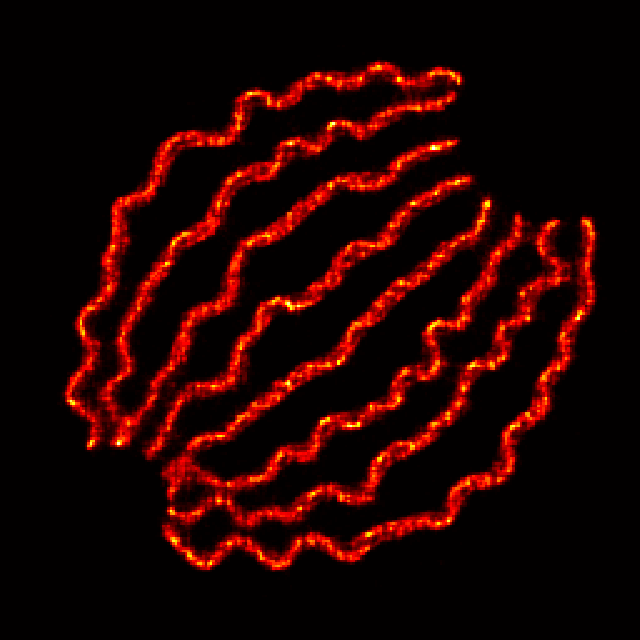} 
        \caption*{$J_4 = 0.3378$}
    \end{subfigure}\\
        \begin{subfigure}{0.3\linewidth}
        \centering
        \includegraphics[width = 0.7\textwidth]{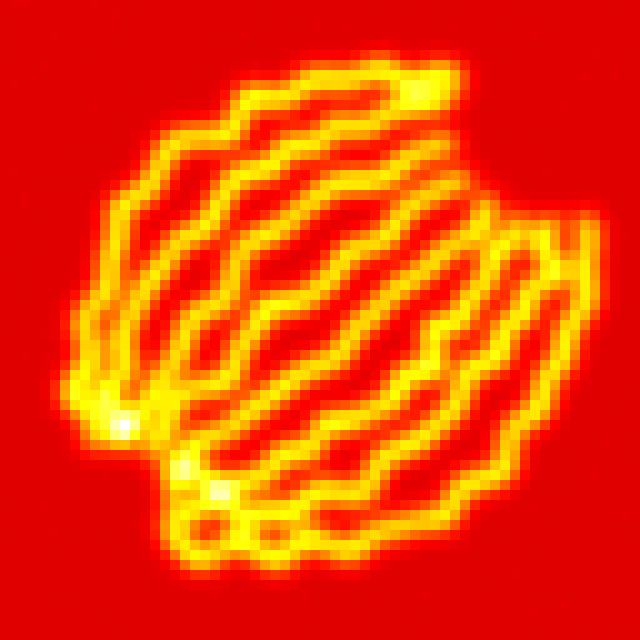} 
        \caption*{}
    \end{subfigure}%
    \begin{subfigure}{0.3\linewidth}
        \centering
        \includegraphics[width = 0.7\textwidth]{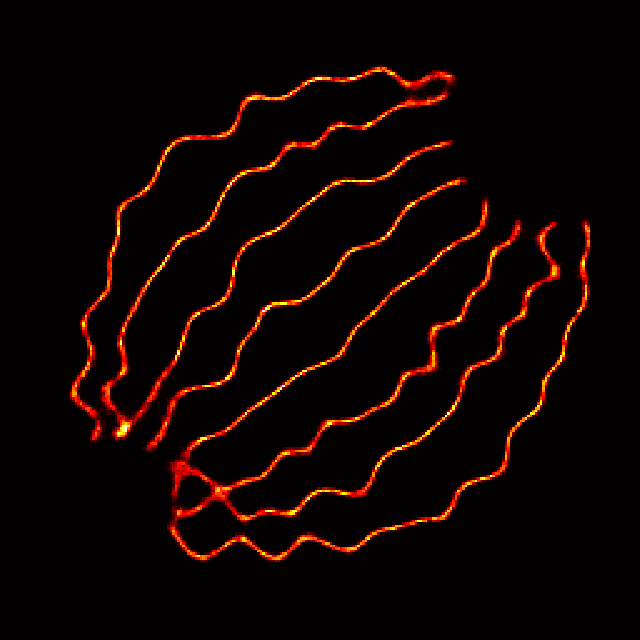} 
        \caption*{$J_4 = 0.6569$}
    \end{subfigure}%
    \begin{subfigure}{0.3\linewidth}
        \centering
        \includegraphics[width = 0.7\textwidth]{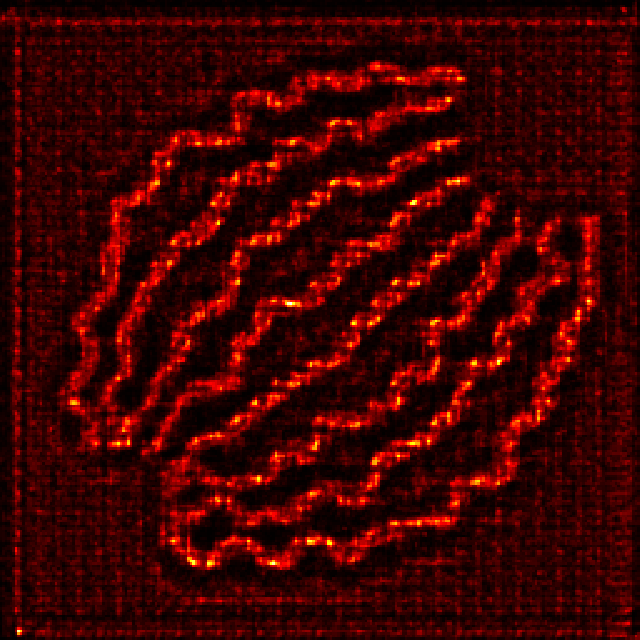} 
        \caption*{\small{$J_4 = 0.1133$}}
    \end{subfigure}\\
    \begin{subfigure}{0.3\linewidth}
        \centering
        \includegraphics[width = 0.7\textwidth]{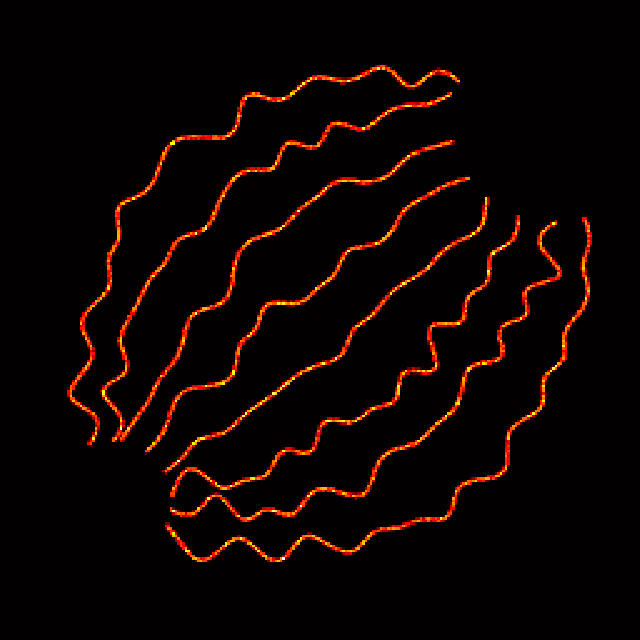}
        \caption{$\overline{g}=\frac{1}{T}\sum_{t=1}^T g_t$}
        \label{fig:exp2_GT}
    \end{subfigure}%
    \caption{Columns from left to right:  low-resolution (widefield) image, reconstruction computed by our approach, DeepSTORM \cite{Nehme_18} reconstruction. Datasets, from top to bottom: Gaussian, Poisson noise with low background and Poisson noise with high background. (a) represents the ground truth stack.
    }
    \label{fig:exp2_recs}
\end{figure}

In the following experiments, we set $c = \delta + 0.01$ and used only $20$ images of the stack for training. The number of inner iterations was set to $K=300$. All the results are reported in Tables \ref{tab:exp2_part1}, \ref{tab:exp2_part2} and \ref{tab:exp2_part3}. Both the widefield-type images 
of the ground truth as well as their the acquired data/reconstructed ones for each dataset are reported in Figure \ref{fig:exp2_GT}/\ref{fig:exp2_recs}, respectively. As a comparison, we report in Figure \ref{fig:exp2_recs} (third column) the reconstructions computed by using DeepSTORM \cite{Nehme_18}, an end-to-end deep-learning based procedure properly trained using the available Google COLAB notebook\footnote{\url{https://colab.research.google.com/github/HenriquesLab/ZeroCostDL4Mic/blob/master/Colab_notebooks/Deep-STORM_2D_ZeroCostDL4Mic.ipynb}} under a physical parameter setting compatible with the datasets observed. Note that for the DeepSTORM reconstructions the PSNR was computed on the final super-resolved image and not frame-wise as the software outputs only the coordinates of the localizations (used to compute the Jaccard index values) and one final intensity image.
We note in this experiment that  no apparent optimal choice between the Gaussian/Poisson data term can be drawn as similar PSNR and Jaccard scores are obtained on the corresponding reconstructions. We argue that this can be due to the biases introduced by the regularization (which make less important the precise modelling of the noise statistics) and by the overall non-convexity of the learning problem. As far as the DeepSTORM reconstructions are concerned, we observe a general tendency to overestimate the reconstruction support with severe artefacts appearing in the background in the case of high background.  

\begin{table}[!h]
    \centering
    \caption{Results in the case of additive white Gaussian noise dataset. For each quality measure, the best score across the models is in bold. The last column reports the initial parameter configuration.}
    \label{tab:exp2_part1}
    \vskip 3mm
    \begin{tabular}{c | c | c | c | c | c | l}
        \toprule
        Model    & $J_0$     & $J_2$     & $J_4$     & Avg. J & PSNR & Initial $\rho$, $\alpha$, $\delta$ \\
        \midrule
        $\mathcal{L}_2$  & 0.1147 & 0.2175 & 0.2194 & 0.1839 &\textbf{ 34.27} & 4, 0.1, - \\
        \midrule
        $\mathcal{L}_1$  & \textbf{0.1309} & \textbf{0.4889} & \textbf{0.5921} & \textbf{0.4040} & 22.45 & 4, 0.1, 5 \\
        \midrule
        DeepSTORM        & 0.0453 & 0.4563 & 0.5688 & 0.3568 & 20.10 & - \\
        \bottomrule
    \end{tabular}
\end{table}

\begin{table}[!h]
    \centering
    \caption{Results in the case of low Poisson noise dataset. For each column, the best score across the models is in bold. The last column reports the initial parameter configuration.} 
    \label{tab:exp2_part2}
    \vskip 3mm
    \begin{tabular}{c c | c | c | c | c | c | l}
        \toprule
        Model   & Fidelity       & $J_0$      & $J_2$     & $J_4$     & Avg. J & PSNR & Initial $\rho$, $\alpha$, $\delta$ \\
        \midrule
        $\mathcal{L}_2$ & KL             & 0.0921  & 0.2633 & 0.3159 & 0.2238 & 33.94 & 50, 5$\cdot 10^{-4}$, - \\
        $\mathcal{L}_2$ & $\|\cdot\|_2$  & 0.1057  & 0.1205 & 0.1211 & 0.1158 & \textbf{34.74} & 4, 0.075, - \\ 
        \midrule
        $\mathcal{L}_1$ & KL             & 0.0907  & 0.4250 & \textbf{0.5924} & 0.3694 & 22.66 & 50, 5$\cdot 10^{-4}$, 0.05\\
        $\mathcal{L}_1$ & $\|\cdot\|_2$  & \textbf{0.1287}  & 0.4494 & 0.5358 & \textbf{0.3713} & 22.52 & 4, 0.1, 5 \\
        \midrule
        DeepSTORM       & - & 0.0381 & \textbf{0.4568} & 0.5890 & 0.3613 & 19.52 & - \\
        \bottomrule
    \end{tabular}
\end{table}

\begin{table}[!h]
    \centering
    \caption{Results in the case of high Poisson noise dataset. For each column, the best score across the models is in bold. The last column reports the initial parameter configuration.} 
    \label{tab:exp2_part3}
    \vskip 3mm
    \begin{tabular}{c c | c | c | c | c | c | l }
        \toprule
        Loss function   & Fidelity       & $J_0$ & $J_2$ & $J_4$ & Avg. J & PSNR & Initial $\rho$, $\alpha$, $\delta$ \\
        \midrule
        $\mathcal{L}_2$ & KL             & 0.0589 & 0.0909 & 0.0987 & 0.0828 & 34.35 & 1, 1, n.a.\\
        $\mathcal{L}_2$ & $\|\cdot\|_2$  & 0.0544 & 0.0778 & 0.0829 & 0.0717 & \textbf{34.36} & 20, 0.05, - \\
        \midrule
        $\mathcal{L}_1$ & KL             & \textbf{0.0982} & 0.3684 & 0.4771 & 0.3146 & 25.31 & 0.5, 0.5, 1 \\
        $\mathcal{L}_1$ & $\|\cdot\|_2$  & 0.0950 & \textbf{0.3902} & \textbf{0.5012} & \textbf{0.3288} & 25.74 & 20, 0.01, 3 \\
        \midrule
        DeepSTORM       & - & 0.0106 & 0.1949 & 0.3185 & 0.1747 & 17.56 & - \\
        \bottomrule
    \end{tabular}
\end{table}

\medskip

We now consider a different ISBI dataset \footnote{ \url{https://srm.epfl.ch/DatasetPage?name=MT1.N1.LD}}, where ground truth images are obtained by discretizing 2D positions of 8731 fluorescent molecules among which approximately 20 molecules are randomly activated at each frame for $436$ ground truth training images. To simulate a long emission (ON) time $\tau_{\text{em}}$, we considered in this experiment a value $\tau_{\text{em}}>\tau_{\text{acq}}$, corresponding to the temporal resolution of the acquisitions. 
When active, a fluorescent molecule is assigned a value in $u=\{0,1,2,3,4\}$ with probability density $p =[0.1,0.25,0.3,0.25,0.1]$ so that for $i\in u$, $p(i)=p_i$. The same Gaussian PSF as in the previous experiment was considered with a medium Poisson noise modelling with a background intensity $b=10$. The binarization was fully learned again with $c = \delta + 0.01$. The number of unrolled iterations of $\cA$ was set to $K = 300$. 

\begin{figure}[!ht]
    \centering
    \begin{tabular}{l l}
         \includegraphics[scale = 0.3]{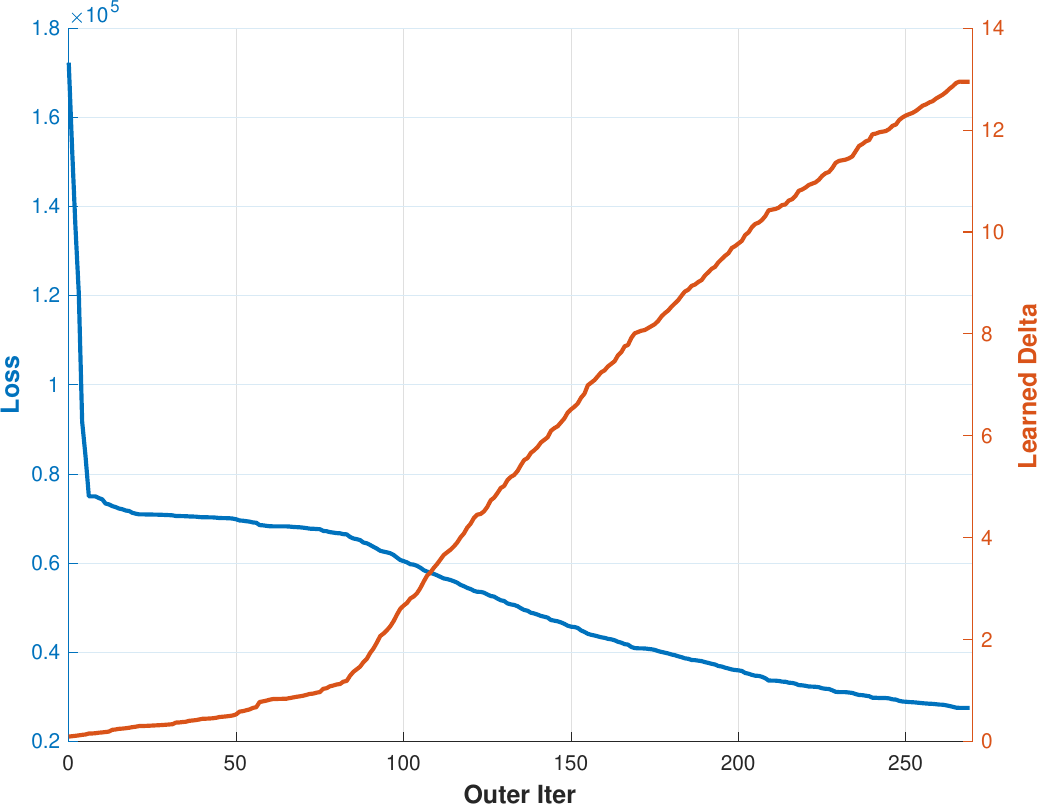} & \includegraphics[scale = 0.3]{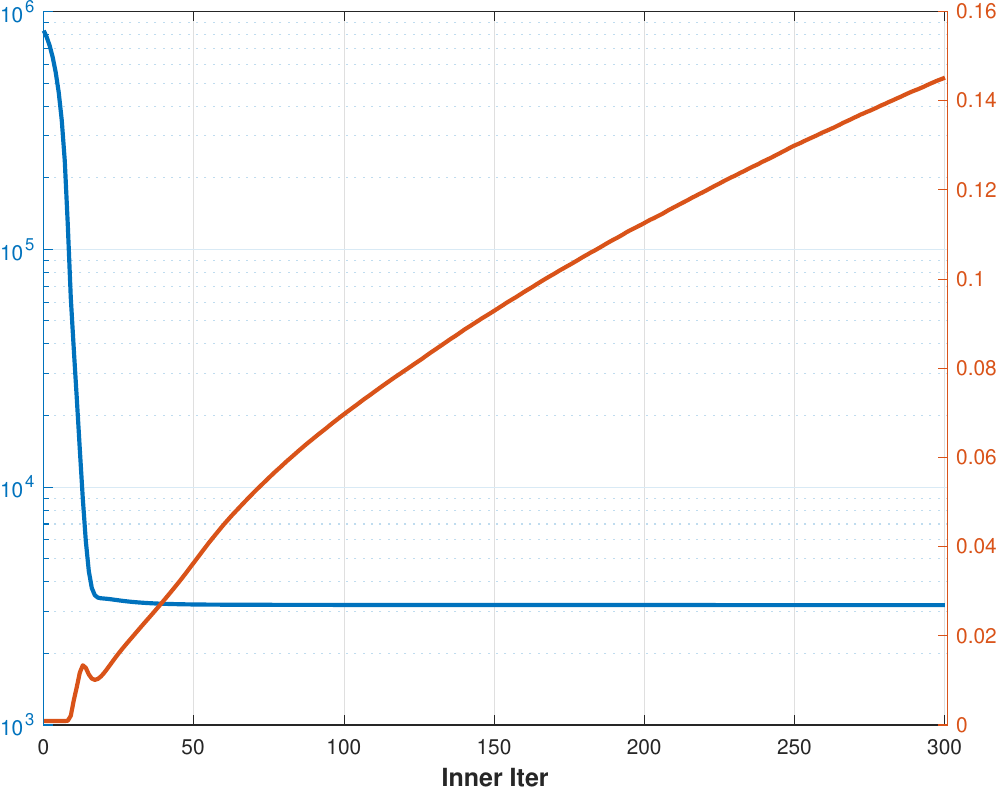} 
    \end{tabular}
    \caption{Left: decrease of the loss functions $\cL_1$ through the outer iterations (blue) and value of the binarization parameter $\delta$ (orange). Right: decrease of the inner energy functional (blue) and increase of Jaccard index (orange), through the inner iterations, computed as the mean across the 416 test samples, without applying the binarization.} 
    \label{fig:exp3_UnrollingPlot}
\end{figure}

\begin{figure}[!ht]
  \centering
  \begin{subfigure}{.24\linewidth}
    \centering
    \includegraphics[width = 0.95\textwidth]{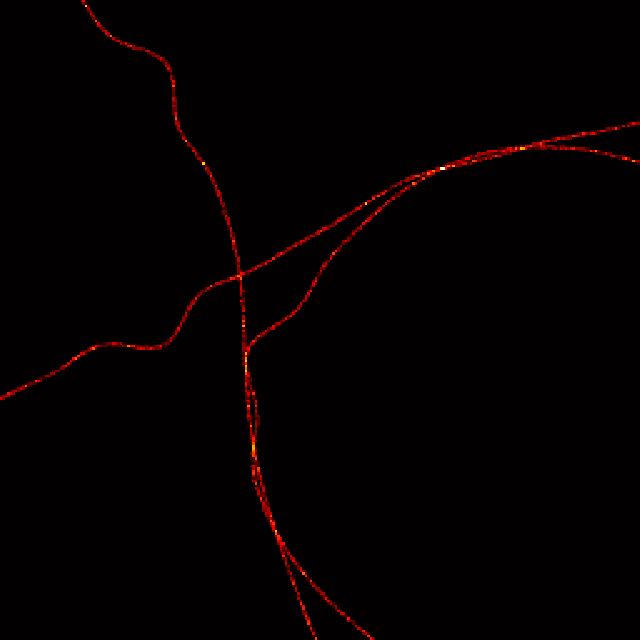}
    \caption{} \label{fig:exp3_GT}
  \end{subfigure}
  \begin{subfigure}{.24\linewidth}
    \centering
    \includegraphics[width = 0.95\textwidth]{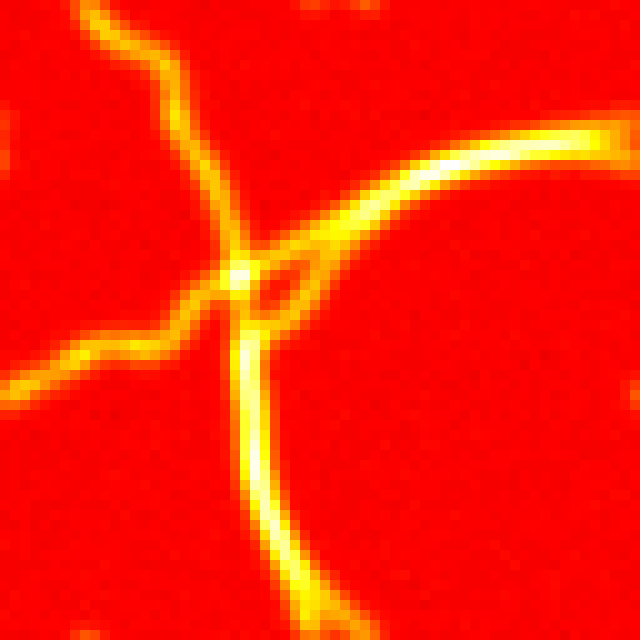}
    \caption{}
  \end{subfigure}
  \begin{subfigure}{.24\linewidth}
    \centering
    \includegraphics[width = 0.95\textwidth]{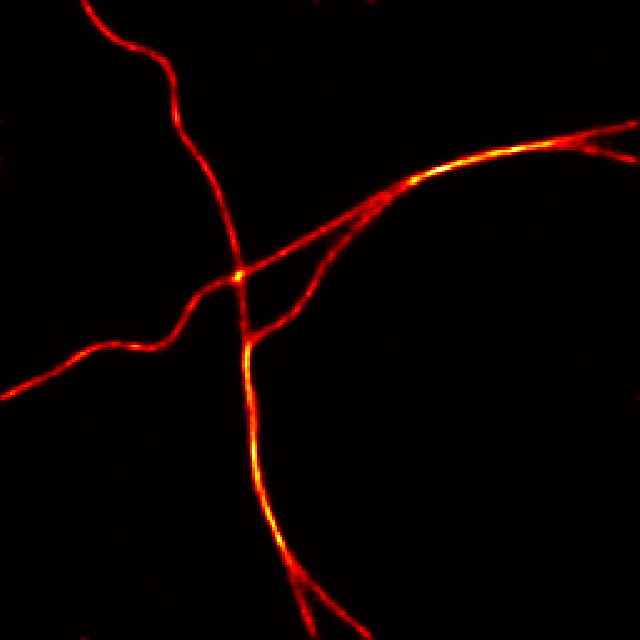}
    \caption{}
  \end{subfigure}
    \begin{subfigure}{.24\linewidth}
    \centering
    \includegraphics[width = 0.95\textwidth]{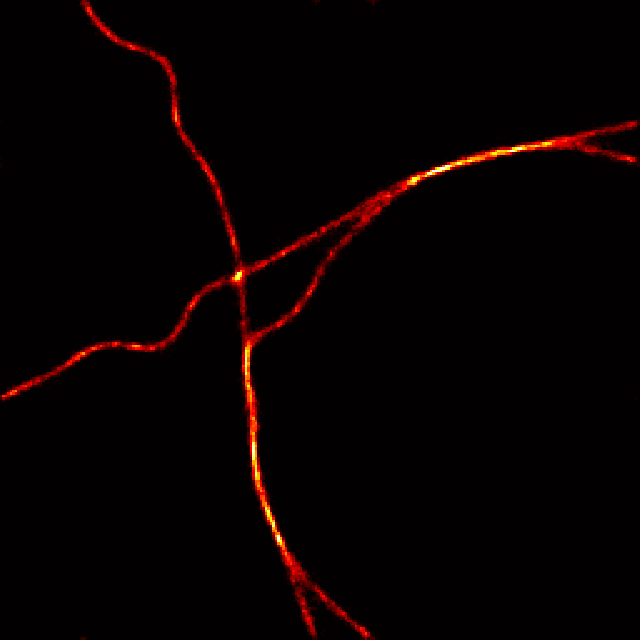}
    \caption{}
  \end{subfigure}
  \caption{(a) Ground truth image. (b) Blurred, noisy and low-resolution measurements (average over time). (c) Reconstruction computed using the $\cL_2$ loss for training. (d) Reconstruction computed using the $\cL_1$ loss for training.} 
    \label{fig:exp3_recs}
\end{figure}

Table \ref{tab:exp3} shows the Jaccard index and PSNR values for this test, while in Figure \ref{fig:exp3_recs} we report the reconstructed images. In this test the use of the KL modelling outperforms a Gaussian fidelity. In Figure \ref{fig:exp3_UnrollingPlot} we report the decay of the loss function \eqref{eq:LossL1} and the value of $\delta$ through the outer iterations. In the other graph, we report the decay of the lower-level energy functional during testing. We observe that the choice $K = 300$ iterations is good enough to guarantee the increasing of the Jaccard index.
\begin{table}[!ht]
    \centering
    \caption{Results for medium Poisson noise dataset. For each column, the best score across the models is in bold. The last column reports the initial parameter configuration.} 
    \label{tab:exp3}
    \vskip 3mm
    \begin{tabular}{c c | c | c | c | c | c | l}
        \toprule
        Loss function   & Fidelity       & $J_0$ & $J_2$ & $J_4$ & Avg. J & PSNR & Initial $\rho$, $\alpha$, $\delta$ \\
        \midrule
        $\mathcal{L}_2$ & KL             & 0.0619 & 0.1285 & 0.1363 & 0.1089 & 39.88 & 1, 0.5, n.a. \\
        $\mathcal{L}_2$ & $\|\cdot\|_2$  & 0.0545 & 0.0812 & 0.0837 & 0.0731 & \textbf{39.97} & 10, 0.01, n.a.\\
        \midrule
        $\mathcal{L}_1$ & KL             & \textbf{0.0776} & \textbf{0.3687} & \textbf{0.4995} & \textbf{0.3153} & 27.39 & 0.1, 0.5, 0.1  \\
        $\mathcal{L}_1$ & $\|\cdot\|_2$  & 0.0661 & 0.2291 & 0.3259 & 0.2070 & 28.00 & 1, 0.05, 0.1 \\
        \bottomrule
    \end{tabular}
\end{table}

\subsection{Semi-blind image  deconvolution in fluctuation-based microscopy.}\label{subsec:Exp4}

We now consider a slightly different image formation model and test the proposed unrolling strategy on an image deconvolution problem based on the study of second-order image fluctuation statistics \cite{SPARCOM,Stergiopoulou2021,stergiopoulou2022,Stergiopoulou_2023_Col0rmeSSVM} with semi-blind estimation of the underlying PSF. Starting from \cite{SOFI}, this approach has proved indeed effective when dealing with biological samples featuring a high density of standard fluorescent molecules which, differently from the ones employed in SMLM applications, are less harmful for the samples observed and thus more suited for \emph{in vivo} analyses. By exploiting the mutual independence of the intensity fluctuations at each pixel and computing second- (or higher-)order statistical information, such approaches reformulate the underlying inverse problem in a covariance domain, which, by the modelling assumptions, enforces a shrinkage of the optical PSF and where sparse regularization can be enforced.

\begin{figure}[!ht]
    \centering
    \includegraphics[scale = 0.35]{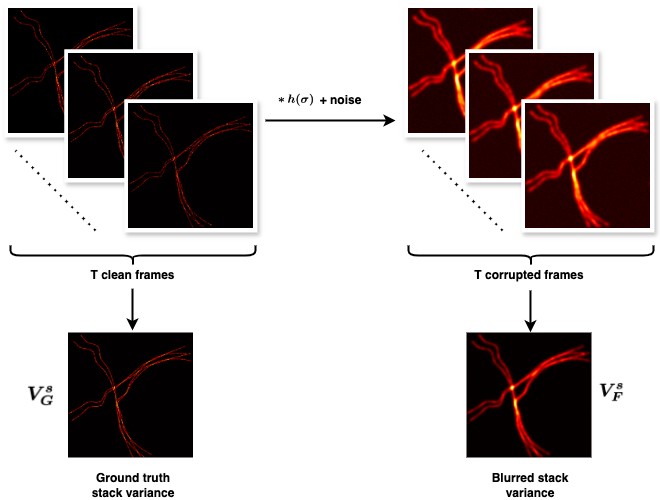}
    \caption{Visualization of the acquisition process for one data sample in the Experiment 4 scenario.}
    \label{fig:DiagrammaExp4}
\end{figure}

The training dataset $\cD = \{(V^s_G,V^s_F)\}_{s = 1}^S$ is generated as follows:
\begin{itemize}
    \item[-] Step 1: for each simulated spatial pattern $s=1,\ldots,S$ a collection of $T = 1000$ fluctuating images is obtained for the molecular structure of interest, $\{g_s^{(t)}\}_{t = 1}^T$. 
    For simulating fluorescence fluctuations over time the SOFI simulator tool \cite{SOFI} was employed.   
    \item[-] Step 2: for each $s=1,\ldots,S$ and each frame  $t=1,\ldots,T$, $g_s^{(t)}$ is blurred with a Gaussian PSF of standard deviation $\varsigma = 3$ and additive white Gaussian noise $n\sim \mathcal{N}(0,\sigma^2 \text{Id})$ with $\sigma = 3$ thus obtaining a  measured image:
    \begin{equation}   \label{eq:invProb_fluct}
    f_s^{(t)} = H(\varsigma)g_s^{(t)} + n.
    \end{equation}
    \item[-] Step 3: for each $s=1,\ldots,S$ the empirical auto-covariances (that is, the variances) $V_G^s$ and  $V_F^s$ of the stack are computed for both the clean and noisy data, respectively, by the formulas: 
    \begin{align*}
        &M_G^s  = \frac{1}{T} \sum_{t = 1}^T g_s^{(t)}, \quad 
        V_G^s  = \frac{1}{T-1} \sum_{t = 1}^T \left(g_s^{(t)} - M_G^s\right)^2; \\
        &M_F^s  = \frac{1}{T} \sum_{t = 1}^T f_s^{(t)}, \quad
        V_F^s  = \frac{1}{T-1} \sum_{t = 1}^T \left(f_s^{(t)} - M_F^s\right)^2 .
    \end{align*}
\end{itemize}
For each spatial pattern $s=1,\ldots,S$ the pair  $(V^s_G,V^s_F)$ thus has as a ground-truth image $V^s_G$ the variance of the fluctuating sample and its corrupted version $V^s_F$  where by \eqref{eq:invProb_fluct}
  \begin{equation}   \label{eq:PSF_shrink}
   V_F^s = H(\varsigma)\odot H(\varsigma)\;V_G^s + V_N,
  \end{equation}
where $\odot$ denotes the point-wise Hadamard product and, by definition, the vector $V_N = \sigma^2\mathbf{e}$ with $\mathbf{e} = (1,\ldots,1)$ denotes the variance of the Gaussian noise considered.
In Figure \ref{fig:DiagrammaExp4} we report a diagram of such modelling.

Note that $V_F^s$ is the second-order SOFI image \cite{SOFI} associated to $\left\{f^{(t)}_s\right\}_{t=1}^T$. By \eqref{eq:PSF_shrink} we notice that the `squaring' of the PSF shrinks its spread by a factor $\sqrt{2}$, thus resulting in a better resolution.
Defining for simplicity $H^2(\varsigma):= H(\varsigma)\odot H(\varsigma) $ the convolution matrix corresponding the kernel $h^2(\varsigma)$. We assume in the following that the Gaussian kernel is defined in terms of an unknown parameter $\sigma$ which we incorporate within the learning procedure, thus considering as lower-level reconstruction functional the model:
\begin{equation}\label{eq:EnergyExp4}
    (V^s_u)^*(\varsigma,\rho)\in \underset{V_u \in \R^{n^2}}{\argmin} \; E_{V_F}(V_u;\theta):=\|H^2(\varsigma) V_u + \sigma^2\mathbf{e} - V^s_F \|_2^2 + \rho \|V_u\|_1 + \iota_{\{V_u \geq 0\}},
\end{equation}
where $\theta$ here includes both the parameters $(\varsigma,\rho)$ and the algorithmic step-sizes $\left\{\alpha_0,\ldots,\alpha_{K-1}\right\}$ over the $K$ unrolled iterations. Note that to backpropagate over $\varsigma$ we have to compute:
\begin{align*}
        \frac{\partial}{\partial \varsigma} \nabla_{V_u} E_{V_F}(V_u;\theta) & = \frac{\partial}{\partial \varsigma} \left( (H^2(\varsigma))^{\mathrm{T}}(H^2(\varsigma) V_u + \sigma^2\mathbf{e} - V^s_F)\right) \\
        & = \frac{\partial}{\partial \varsigma}(H^2(\varsigma))^\mathrm{T}(H^2(\varsigma) V_u + \sigma^2\mathbf{e} - V^s_F) + (H^2(\varsigma))^\mathrm{T} \left(\frac{\partial}{\partial \varsigma}\;H^2(\varsigma) V_u\right),
\end{align*}
Further computations show that
\[
\frac{\partial}{\partial \varsigma}\;H^2(\varsigma) V_u = \left( \frac{\partial}{\partial \varsigma}H^2(\varsigma)\right)V_u = \frac{\partial}{\partial \varsigma} h^2(\varsigma) \ast V_u,
\]
with
\[
        \frac{\partial}{\partial \sigma}h^2(\varsigma)[x,y]  = \frac{e^{-\frac{\left(x - \frac{n+1}{2}\right)^2 + \left(y - \frac{n+1}{2}\right)^2}{\varsigma^2}}}{2\pi^2\varsigma^5}
        \left(\frac{1}{\varsigma^2}\left(\left(x - \frac{n+1}{2}\right)^2 + \left(y - \frac{n+1}{2}\right)^2\right) - 2\right).
\]

For constructing the dataset $\cD$ we considered $S=30$ different spatial patterns. Each pattern is obtained by super-position/rotation and translation of filament structures similar to the ones in Figure \ref{fig:exp3_GT}.
The dataset was split into a training set with 20 elements and a test set with 10 images. The same experimental setup was used for the lower-level solver, with a number of $K=300$ unrolled iterations. The $\cL_1$ penalty \eqref{eq:LossL1} was used for training. The binarization function was fully learned with $c = \delta + 0.01$. The initial values for the parameters of this experiment were: $\rho = 10^{-5}$, $\delta = 25$, $\alpha = 1000$ and $\varsigma = 5$.

\begin{figure}[!h]
  \centering
  \begin{subfigure}{.24\linewidth}
    \centering
    \includegraphics[width = 0.95\linewidth]{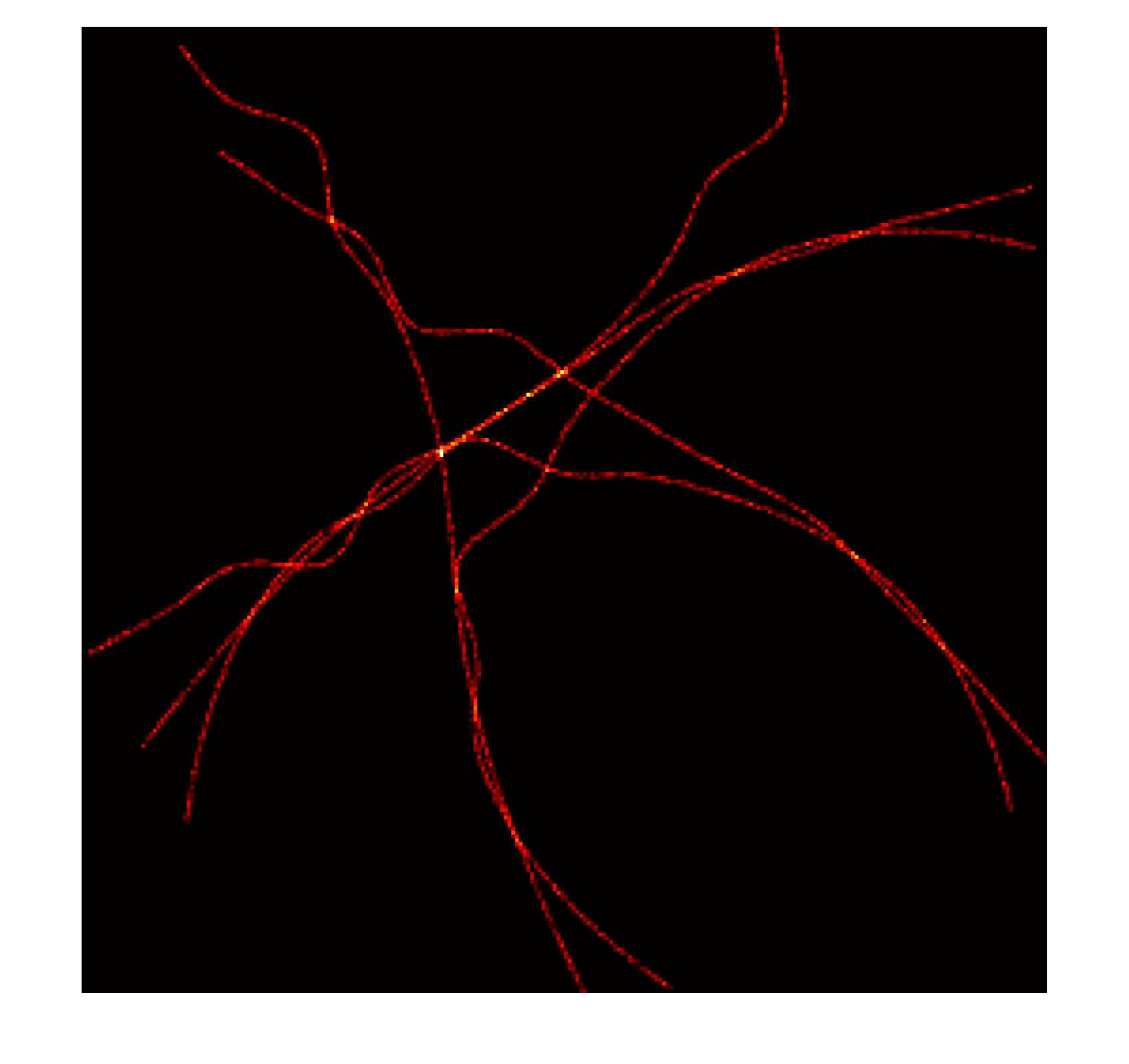}
    \caption{}
  \end{subfigure}%
  \begin{subfigure}{.24\linewidth}
    \centering
    \includegraphics[width = 0.95\linewidth]{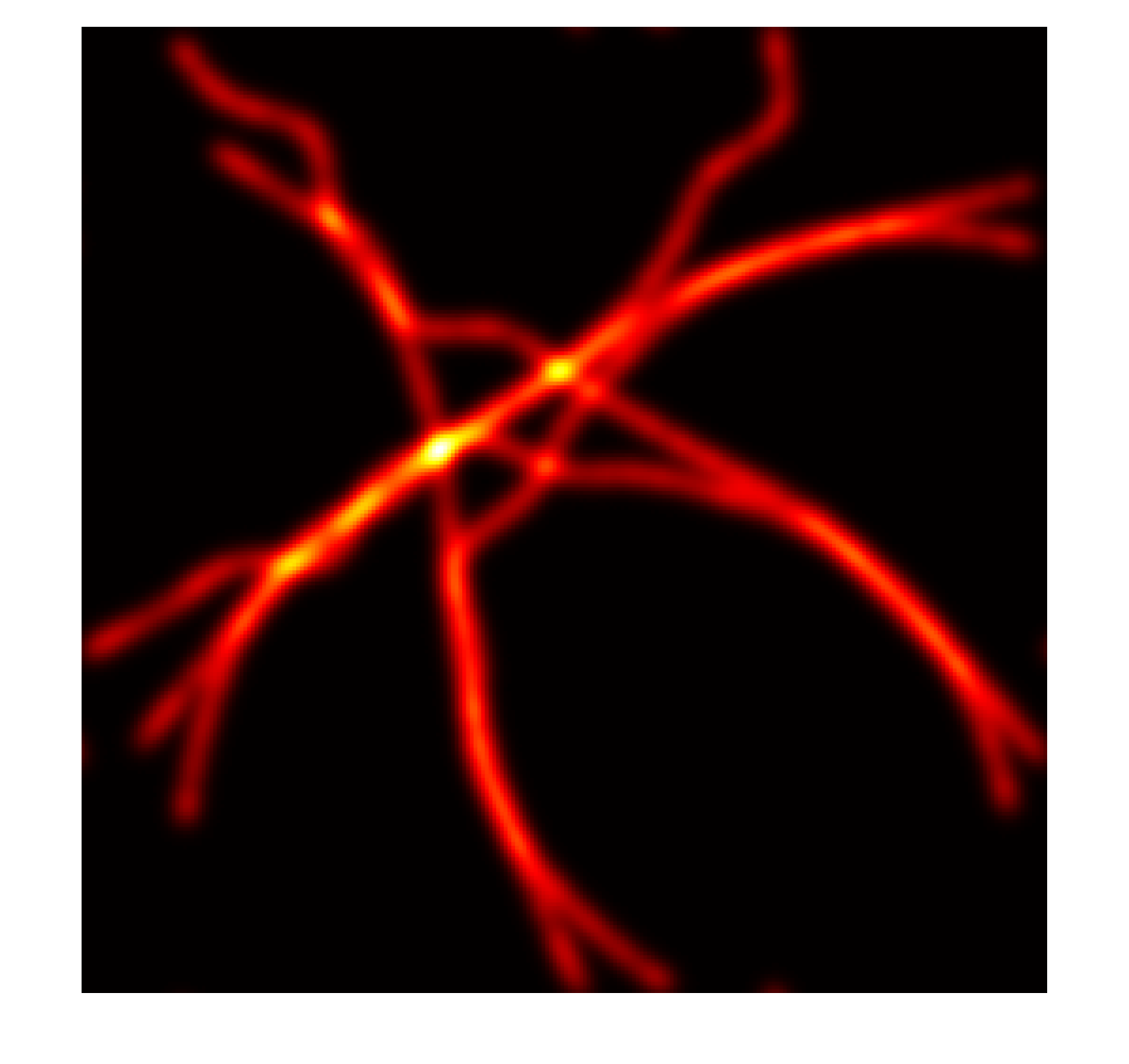}
    \caption{}
  \end{subfigure}%
  \begin{subfigure}{.24\linewidth}
    \centering
    \includegraphics[width = 0.95\linewidth]{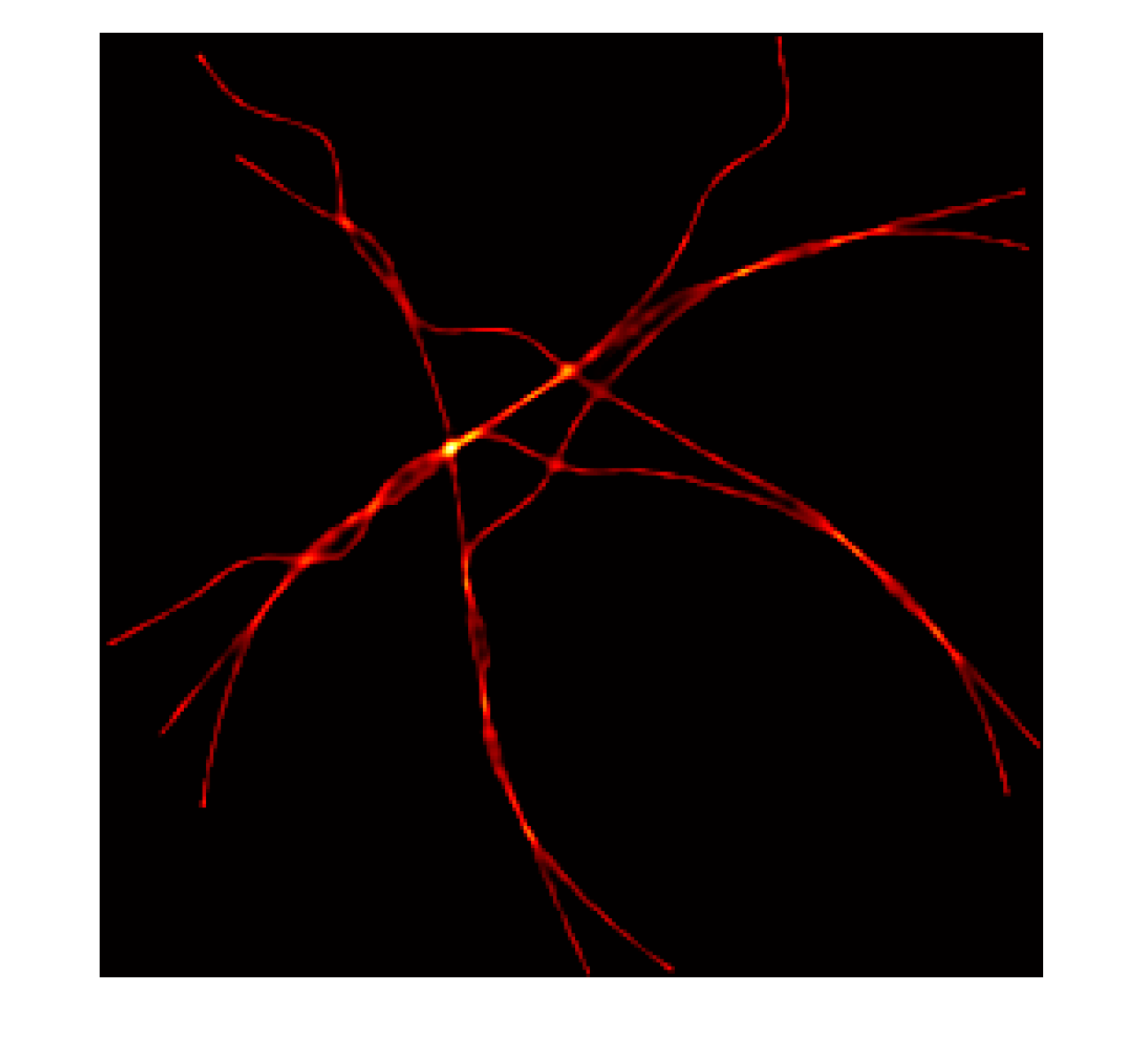}
    \caption{}
  \end{subfigure}
    \begin{subfigure}{.24\linewidth}
    \centering
    \includegraphics[width = 0.95\linewidth]{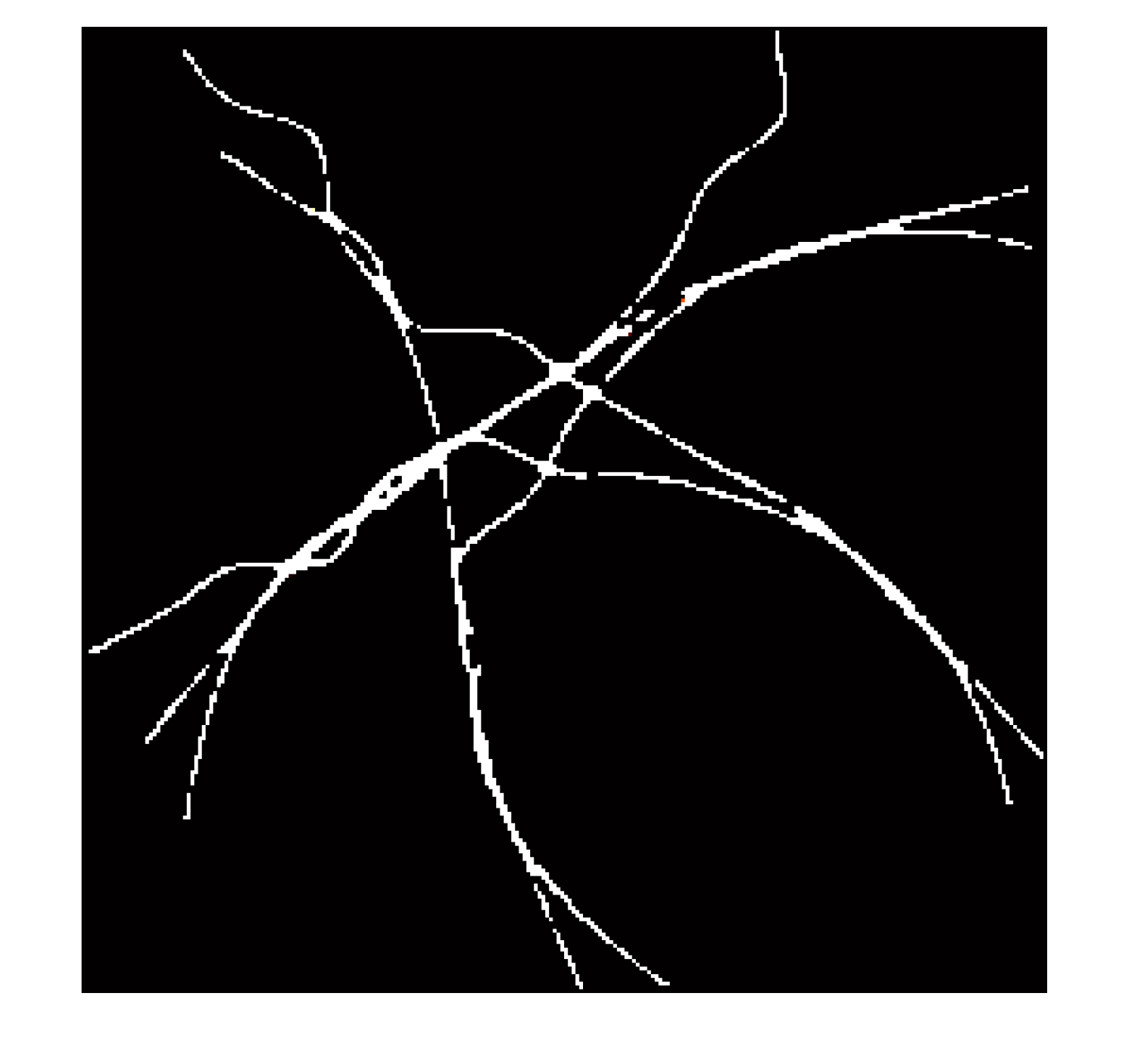}
    \caption{}
  \end{subfigure}
  \caption{(a) Ground truth image $V_{G}^s$. (b) Noisy variance image $V_{F}^s$. (c) Computed reconstruction. (d) Computed reconstruction post-binarization.}
    \label{fig:exp4_recs}
\end{figure}

The average Jaccard index computed over the whole test set is reported in Table \ref{tab:exp4}. Values of the Jaccard index computed before the binarization is applied are also reported to show the beneficial effect of introducing such function within the training phase. The learned value for $\varsigma$ was $\varsigma^{\ast} = 2.87$. In Figure \ref{fig:exp4_recs} we illustrate one reconstruction computed from the test dataset corresponding to the optimal parameters learned. Note that the deconvolution performance appears quite accurate. Moreover, despite the slight underestimation of the PSF, no ringing effect is observed, probably due to the sparsity induced by the use of the $\ell_1$ regularization with the non-negativity constraint.

\begin{table}[!h]
    \centering
    \caption{Localization assessment for the semi-blind deconvolution problem, with/without binarization. 
    } 
    \label{tab:exp4}
    \vskip 3mm
    \begin{tabular}{c| c | c | c | c}
        \toprule
        Binarization  & $J_0$     & $J_2$     & $J_4$     & Avg. J \\
        \midrule
            No        & 0.7096 & 0.7164 & 0.7208 & 0.7156 \\
        \midrule
            Yes       & 0.7922 & 0.8362 & 0.8582 & 0.8289 \\
        \bottomrule
    \end{tabular}
\end{table}

\begin{figure}[!h]
    \centering
    \begin{tabular}{l l}
         \includegraphics[scale = 0.3]{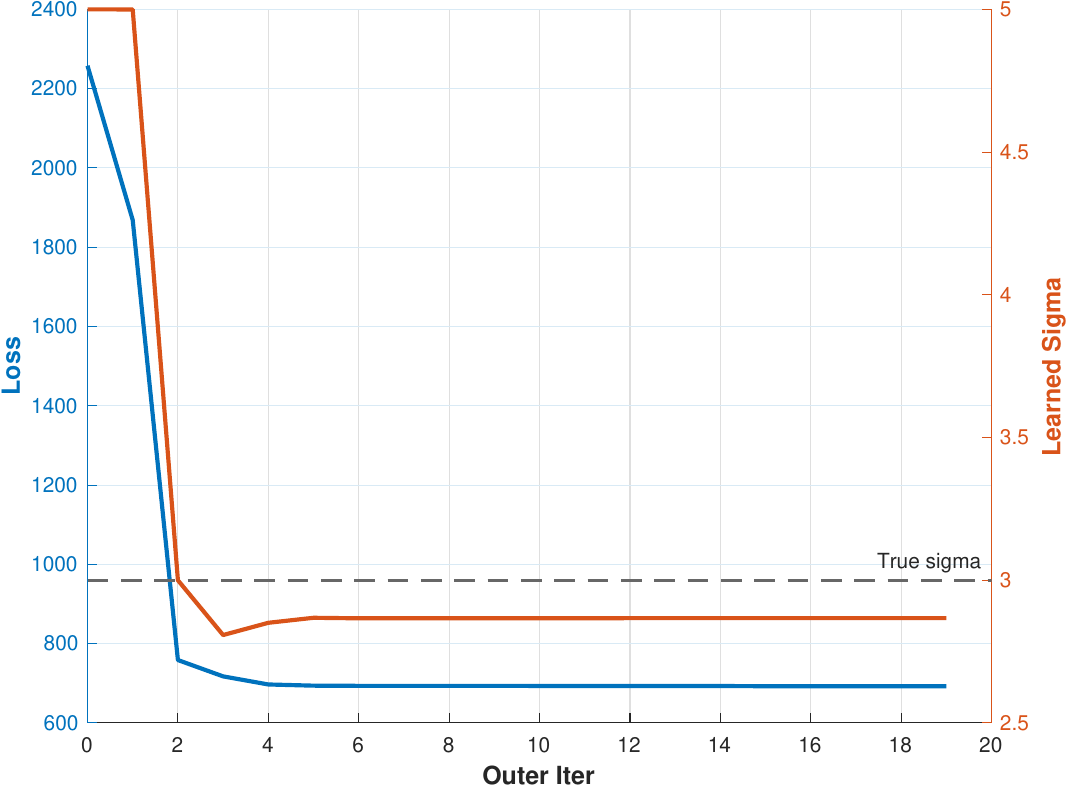} & \includegraphics[scale = 0.3]{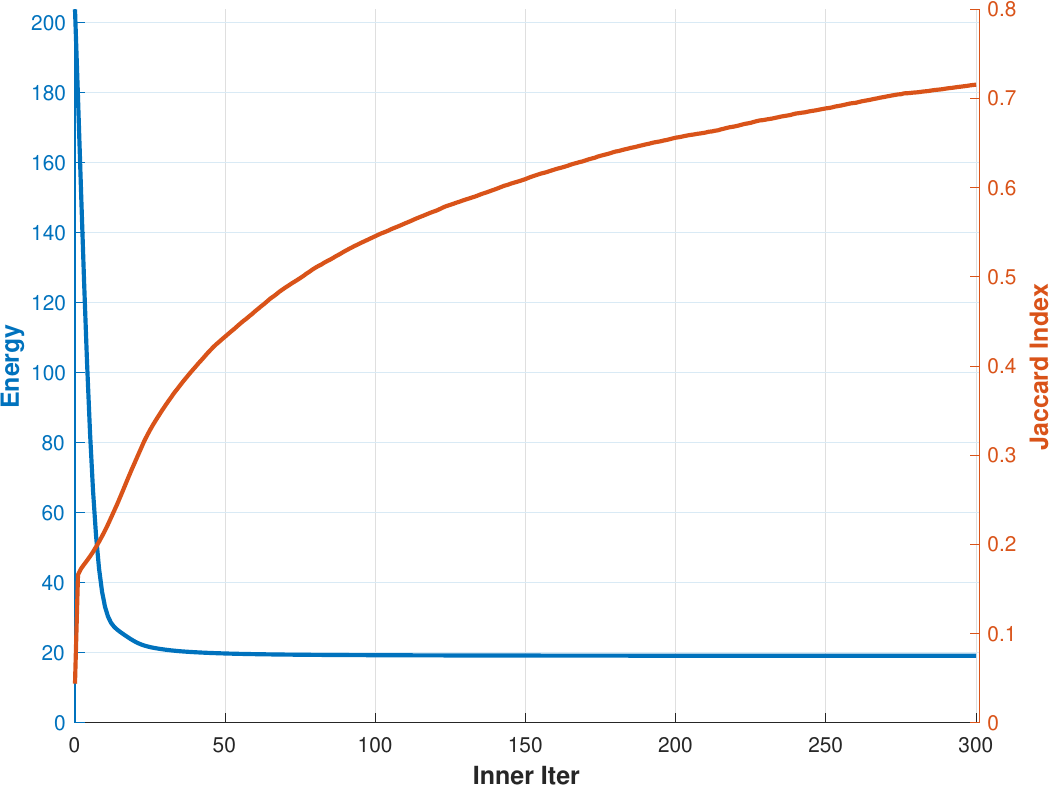} 
    \end{tabular}
    \caption{Left: decrease of the $\mathcal{L}_1$ loss function through the outer iterations (blue) and evolution of the learned $\varsigma$ (orange). Right image: decrease of the inner energy functional \eqref{eq:EnergyExp4} (blue) and increase of Jaccard index (orange) along inner iterations, computed as the mean across the 10 test samples.
    }
    \label{fig:exp4_UnrollingPlot}
\end{figure}

\section{Conclusions}

We proposed an unfolded accelerated projected gradient descent scheme to solve both reconstruction and localization inverse problems often arising in the field of fluorescence microscopy.

By considering  $\ell_1$-regularized criteria featuring suitable (Gaussian, Poisson) data terms, the proposed strategy learn optimal regularization and algorithmic parameters by means of two different loss functions assessing reconstruction and localization performance by suitable comparisons with given ground truth data.

The framework proposed adapts to both SMLM and fluctuation-based approaches as showed by several simulated and realisitic numerical results, even in the presence of partially known physical models which can be learned within the training procedure.

Further work should address and compare the performance of unrolled schemes when used either in a physics-driven or in an end-to-end fashion as well as potentially combined with more effective, even neural network-based, sparsity promoted regularizations, such in the case of Plug and Play regularization procedures.


\section*{Acknowledgments}
The authors acknowledge Vasiliki Stergiopoulou for kindly providing the datasets and the specifics for the tests run in Section \ref{subsec:Exp4}.
Silvia Bonettini, Danilo Pezzi and Marco Prato are all members of the INdAM research group GNCS, which is kindly acknowledged.

\section*{Declarations}

\noindent
\paragraph{Funding.} SB, MP, DP are partially supported by the PRIN projects 2022ANC8HL and 20225STXSB, under the National Recovery and Resilience Plan (NRRP) funded by the
European Union – NextGenerationEU. LC acknolwedges the support received by the ANR JCJC TASKABILE project (ANR-22-CE48-0010) and by the ANR PRC MICROBLIND (ANR-21-CE48-0008).

\paragraph{Competing interests.} The authors have no competing interests to declare that are relevant to the content of this article.

\paragraph{Ethics approval.} The authors declare that research ethics approval was not required for this study.

\paragraph{Consent.} The authors declare that informed consents were not required for this study.

\paragraph{Data availability statement.} 
The images of Section \ref{subsec:Exp2} can be downloaded from the websites \url{https://srm.epfl.ch/srm/index.php} and \url{https://srm.epfl.ch/DatasetPage?name=MT1.N1.LD}. All the synthetic datasets used in the making of this paper can be found at \url{https://github.com/DaniloPezzi/AGD-unfolding-for-microscopy-imaging}.

\paragraph{Conflict of interests.} All authors certify that they have no affiliations with or involvement in any organization or entity with any financial interest or non-financial interest in the subject matter or materials discussed in this manuscript.

\paragraph{Authors' contributions.} All authors contributed to the study conception and design of the paper. The unrolled scheme has been created by S. Bonettini, D. Pezzi and M. Prato. The application to the SMLM problem has been handled by L. Calatroni. Code implementation, data collection and analysis were performed by L. Calatroni and D. Pezzi. The first draft of the manuscript was written by D. Pezzi, and all authors commented on previous versions of the manuscript. All authors read and approved the final manuscript.

\appendix
\section{Non-differentiability of CEL0 w.r.t.~the regularization parameter} \label{app:A}

We mentioned in Section \ref{sec:modelSparse} that the differentiability w.r.t.~the regularization parameter denoted by $\rho$ in \eqref{eq:E} is required to apply the proposed unrolling strategy. Here we show such property does not hold for the CEL0 regularizer continuous and non-convex regularization functional $R_{\texttt{CEL0}}$ studied in \cite{Soubies_2015_CEL0,Soubies2017} which, given its performance in localization microscopy problems (see, e.g., \cite{Gazagnes2017,LazzarettiISBI2021}) would be a natural choice for localization purposes.

Such functional takes the form
\begin{equation}
    R_{\texttt{CEL0}}(u) = \sum_{p=1}^n \phi(u_p;\rho,\|a_p\|_2),
\end{equation}
where $\|a_p\|_2$ denotes the norm of the $p$--th column of the linear modeling operator $A$ and $\rho$ is the regularization parameter whose dependence here is highly non-linear. Enforcing non-negativity one has 
\begin{equation}
    \phi(s;\rho,a) =
    \begin{cases}
         \rho - \frac{a^2}{2}\left(s - \frac{\sqrt{2\rho}}{a}\right)^2, & 0\leq s < \frac{\sqrt{2\rho}}{a} \\
         \rho,                                                             & s \geq \frac{\sqrt{2\rho}}{a}
    \end{cases}
\end{equation}
showing that the function $\phi(s;\rho,a)$ is differentiable w.r.t. $s$ with derivative
\begin{equation}
    \phi'(s;\rho,a) = 
    \begin{cases}
         -a^2\left(s - \frac{\sqrt{2\rho}}{a}\right), & 0\leq 
 s < \frac{\sqrt{2\rho}}{a} \ \left(\text{equiv. } \rho > \frac{a^2s^2}{2}\right)\\
         0,                                           & s \geq \frac{\sqrt{2\rho}}{a} \ \left(\text{equiv. } \rho \leq \frac{a^2s^2}{2}\right)
    \end{cases}.
\end{equation}
Since
\begin{equation}
    \lim_{\epsilon \downarrow 0^+}\frac{\phi'(s;\frac{a^2s^2}{2}+\epsilon,a) - \phi'(s;\frac{a^2s^2}{2},a)}{\epsilon} = \lim_{\epsilon \downarrow 0^+} -\frac{a^2}{\epsilon}\left(s - \frac{\sqrt{a^2s^2+2\epsilon}}{a}\right) = \frac{1}{s} \neq 0,
\end{equation}
the function $\phi'(s;\rho,a)$ is not differentiable w.r.t. $\rho$.


\printbibliography

@preamble{ " \newcommand{\noop}[1]{} " }

@article{Helsinki,
	author = {Bonettini, Silvia and Franchini, Giorgia and Pezzi, Danilo and Prato, Marco},
	doi = {10.3934/ipi.2022055},
	journal = {Inverse Probl. Imaging},
	number = {5},
	pages = {925--950},
	title = {Explainable Bilevel Optimization: An Application to the {H}elsinki Deblur Challenge},
	volume = {17},
	year = {2023}
	}

@article{Dickson1997,
	Author = {Dickson, R. M. and Cubitt, A. B. and Tsien, R. Y. and Moerner, W. E.},
	Journal = {Nature},
	Number = {6640},
	Pages = {355--358},
	Title = {On/off blinking and switching behaviour of single molecules of green fluorescent protein},
	Volume = {388},
	Year = {1997},
        DOI = {10.1038/41048}
}

@article{Hell94,
	author = {S. W. Hell and J. Wichmann},
	journal = {Opt. Lett.},
	number = {11},
	pages = {780--782},
	publisher = {Optica Publishing Group},
	title = {Breaking the diffraction resolution limit by stimulated emission: stimulated-emission-depletion fluorescence microscopy},
	volume = {19},
	year = {1994},
        doi = {10.1364/OL.19.000780}
}

@article{Betzig2006,
	author = {E. Betzig  and G. H. Patterson  and R. Sougrat  and O. W. Lindwasser  and S. Olenych  and J. S. Bonifacino  and M. W. Davidson  and J. Lippincott-Schwartz  and H. F. Hess },
	title = {Imaging Intracellular Fluorescent Proteins at Nanometer Resolution},
	journal = {Science},
	volume = {313},
	number = {5793},
	pages = {1642--1645},
	year = {2006},
        doi = {10.1126/science.1127344}
}

@article {iPALM,
	author = {Shtengel, G. and Galbraith, J.A. and Galbraith, C.G. and Lippincott-Schwartz, J. and Gillette, J. M. and Manley, S. and Sougrat, R. and Waterman, C. M. and Kanchanawong, P. and Davidson, M. W. and Fetter, R. D. and Hess, H. F.},
	title = {Interferometric fluorescent super-resolution microscopy resolves {3D} cellular ultrastructure},
	volume = {106},
	number = {9},
	pages = {3125--3130},
	year = {2009},
	doi = {10.1073/pnas.0813131106},
	journal = {Proc. Natl. Acad. Sci. U. S. A.}
}

@article{3D_STORM,
	author = {B. Huang  and W. Wang  and M. Bates  and X. Zhuang },
	title = {Three-Dimensional Super-Resolution Imaging by Stochastic Optical Reconstruction Microscopy},
	journal = {Science},
	volume = {319},
	number = {5864},
	pages = {810-813},
	year = {2008},
	doi = {10.1126/science.1153529},
	URL = {https://www.science.org/doi/abs/10.1126/science.1153529}}

@article{Lelek2021,
	Author = {Lelek, Micka{\"e}l and Gyparaki, Melina T. and Beliu, Gerti and Schueder, Florian and Griffi{\'e}, Juliette and Manley, Suliana and Jungmann, Ralf and Sauer, Markus and Lakadamyali, Melike and Zimmer, Christophe},
	Journal = {Nat. Rev. Method. Prim.},
	Number = {1},
	Pages = {39},
	Title = {Single-molecule localization microscopy},
	Volume = {1},
	Year = {2021},
        DOI = {10.1038/s43586-021-00038-x}
}

@article{Aritake2021,
title = {Fast and robust multiplane single-molecule localization microscopy using a deep neural network},
journal = {Neurocomputing},
volume = {451},
pages = {279-289},
year = {2021},
issn = {0925--2312},
author = {Toshimitsu Aritake and Hideitsu Hino and Shigeyuki Namiki and Daisuke Asanuma and Kenzo Hirose and Noboru Murata},
doi = {10.1016/j.neucom.2021.04.050}
}

@INPROCEEDINGS{Gazagnes2017,
	author={S. {Gazagnes} and E. {Soubies} and L. {Blanc-F\'eraud}},
	booktitle={2017 IEEE 14th International Symposium on Biomedical Imaging}, 
	title={High-density molecule localization for super-resolution microscopy using {CEL0}-based sparse approximation}, 
	year={2017},
	volume={},
        pages={28--31},
        doi={10.1109/ISBI.2017.7950460}
}

@INPROCEEDINGS{LazzarettiISBI2021,
	author={Lazzaretti, M. and  L. Calatroni and Estatico, C.},
	booktitle={2021 IEEE 18th International Symposium on Biomedical Imaging}, 
	title={Weighted-{CEL0} Sparse Regularisation For Molecule Localisation In Super-Resolution Microscopy With {P}oisson Data}, 
	year={2021},
	pages={1751--1754},
        doi = {10.1109/ISBI48211.2021.9434014}
}

@INPROCEEDINGS{Stergiopoulou2021,
	author={Stergiopoulou, V.  and de Morais Goulart, J. H. and Schaub, S. and  L. Calatroni  and Blanc-F\'{e}raud, L.},
	booktitle={2021 IEEE 18th International Symposium on Biomedical Imaging},
	title={{COL0RME}: Covariance-Based $\ell_0$ Super-Resolution Microscopy with Intensity Estimation}, 
	year={2021},
	pages={349--352},
	doi={10.1109/ISBI48211.2021.9433976}}

@article{stergiopoulou2022, 
	title={{COL0RME}: Super-resolution microscopy based on sparse blinking/fluctuating fluorophore localization and intensity estimation}, 
	volume={2}, 
	DOI={10.1017/S2633903X22000010}, 
	journal={Biolog. Imaging}, 
	publisher={Cambridge University Press}, 
	author={Stergiopoulou, V. and Calatroni, L. and de Morais Goulart, H. and Schaub, S. and Blanc-F\'{e}raud, L.}, 
	year={2022}, 
	pages={e1}
}

@article{SPARCOM,
	author = {Solomon, O. and Eldar, Y. C. and Mutzafi, M. and Segev, M.},
	title = {{SPARCOM}: Sparsity Based Super-resolution Correlation Microscopy},
	journal = {SIAM J. Imaging Sci.},
	volume ={12},
	issue ={1},
	pages = {392--419},
	year = {2019},
        doi = {10.1137/18M1174921}
}

@article{srrf,
	author = {Gustafsson, N. and Culley, S. and Ashdown, G. and Owen, D. M. and Pereira, P. M. and Henriques, R.},
	journal = {Nat. Commun.},
        volume = {7},
        number = {1},
	pages = {12471},
	title = {Fast live-cell conventional fluorophore nanoscopy with {ImageJ} through super-resolution radial fluctuations},
	year = {2016},
        doi = {10.1038/ncomms12471}
}

@article {SOFI,
author = {Dertinger, T. and Colyer, R. and Iyer, G. and Weiss, S. and Enderlein, J.},
title = {Fast, background-free, {3D} super-resolution optical fluctuation imaging ({SOFI})},
year = {2009},
doi = {10.1073/pnas.0907866106},
volume = {106},
number = {52},
journal = {Proc. Natl. Acad. Sci. U. S. A.},
pages = {22287--22292}
}

@article{Bonettini_2009_SGP,
	author = {S Bonettini and Zanella, Riccardo and Zanni, Luca},
	doi = {10.1088/0266-5611/25/1/015002},
	journal = {Inverse Probl.},
	number = {1},
	pages = {015002},
	title = {A scaled gradient projection method for constrained image deblurring},
	volume = {25},
	year = {2008}
}

@misc{ghadimi2018approximation,
      title={Approximation methods for bilevel programming}, 
      author={Saeed Ghadimi and Mengdi Wang},
      year={2018},
      eprint={1802.02246},
      archivePrefix={arXiv},
      primaryClass={math.OC},
      note={arXiv:1802.02246},
      doi = {10.48550/arxiv.1802.02246}
}

@ARTICLE{Harmany2012,
  author={Harmany, Zachary T. and Marcia, Roummel F. and Willett, Rebecca M.},
  journal={IEEE Transactions on Image Processing}, 
  title={This is SPIRAL-TAP: Sparse Poisson Intensity Reconstruction ALgorithms—Theory and Practice}, 
  year={2012},
  volume={21},
  number={3},
  pages={1084-1096},
  doi={10.1109/TIP.2011.2168410}}

@inproceedings{NEURIPS2022_aa84ec1a,
	author = {Dagr\'{e}ou, Mathieu and Ablin, Pierre and Vaiter, Samuel and Moreau, Thomas},
	booktitle = {Advances in Neural Information Processing Systems},
	pages = {26698--26710},
	publisher = {Curran Associates, Inc.},
	title = {A framework for bilevel optimization that enables  stochastic and global variance reduction algorithms},
	volume = {35},
	year = {2022},
    address = {New York}
}

@article{Chambolle_2015_FISTA_Iterates,
	author = {Chambolle, A. and Dossal, Ch.},
	doi = {10.1007/s10957-015-0746-4},
	journal = {J. Optim. Theory Appl.},
	number = {3},
	pages = {968--982},
	title = {On the Convergence of the Iterates of the ``{F}ast {I}terative {S}hrinkage/{T}hresholding {A}lgorithm''},
	volume = {166},
	year = {2015}
	}

@article{Koulouri_2021_AdaptiveSuperRes,
	author = {Koulouri, Alexandra and Heins, Pia and Burger, Martin},
	doi = {10.1109/TSP.2020.3037373},
	journal = {IEEE Trans. Signal Process.},
	pages = {165--178},
	title = {Adaptive Superresolution in Deconvolution of Sparse Peaks},
	volume = {69},
	year = {2021}
	}

@ARTICLE{Bertero-etal-2009,
	author={Bertero, M. and Boccacci, P. and Desider\`{a}, G. and Vicidomini, G.},
	title={Image deblurring with {P}oisson data: From cells to galaxies},
	journal={Inverse Probl.},
	year={2009},
	volume={25},
	number={12},
	pages={123006},
        doi={10.1088/0266-5611/25/12/123006}
}

@article{Toader2022,
	author = {Toader, B. and Boulanger, J. and Korolev, Y. and Lenz, M. O. and Manton, J. and Sch{\"o}nlieb, C.-B. and Mure{\c s}an, L.},
	doi = {10.1007/s10851-022-01100-3},
	journal = {J. Math. Imaging Vis.},
	title = {Image Reconstruction in Light-Sheet Microscopy: Spatially Varying Deconvolution and Mixed Noise},
	year = {2022},
        volume = {64},
        pages = {968--992}
}

@article{CalatroniSIAM2017,
	author = { L. Calatroni  and De Los Reyes, J. C. and Sch\"{o}nlieb, C.-B.},
	title = {Infimal Convolution of Data Discrepancies for Mixed Noise Removal},
	journal = {SIAM J. Imaging Sci.},
	volume = {10},
	number = {3},
	year = {2017},
        pages = {1196-1233},
        doi = {10.1137/16M1101684}
}

@article{ROF,
	title = "Nonlinear total variation based noise removal algorithms",
	journal = "Physica D: Nonlinear Phen.",
	volume = "60",
	number = "1",
	pages = "259--268",
	year = "1992",
	issn = "0167-2789",
	doi = "10.1016/0167-2789(92)90242-F",
	author = "L. I. Rudin and S. Osher and E. Fatemi"
}

@ARTICLE{Jezierska2014,
	author={A. {Jezierska} and C. {Chaux} and J. {Pesquet} and H. {Talbot} and G. {Engler}},
	journal={IEEE Trans. on Sign. Process.}, 
	title={An {EM} Approach for Time-Variant {P}oisson-{G}aussian Model Parameter Estimation}, 
	year={2014},
	volume={62},
	number={1},
        pages = {17--30},
        doi = {10.1109/TSP.2013.2283839}
}

@inproceedings{GregorLeCun2010,
author = {Gregor, Karol and LeCun, Yann},
title = {Learning fast approximations of sparse coding},
year = {2010},
isbn = {9781605589077},
publisher = {Omnipress},
booktitle = {Proceedings of the 27th International Conference on International Conference on Machine Learning},
pages = {399–-406},
address = {Madison},
}

@ARTICLE{Adler2018,
  author={Adler, Jonas and \"{O}ktem, Ozan},
  journal={IEEE Trans. Med. Imaging}, 
  title={Learned Primal-Dual Reconstruction}, 
  year={2018},
  volume={37},
  number={6},
  pages={1322--1332},
  doi={10.1109/TMI.2018.2799231}
}

@article{Bertocchi_2020,
year = {2020},
publisher = {IOP Publishing},
volume = {36},
number = {3},
pages = {034005},
author = {C Bertocchi and E Chouzenoux and M-C Corbineau and J-C Pesquet and M Prato},
title = {Deep unfolding of a proximal interior point method for image restoration},
journal = {Inverse Probl.},
doi = {10.1088/1361-6420/ab460a}
}

@misc{Brauer2022,
	author = {Brauer, C. and Breustedt, N. and de Wolff, T. and Lorenz, D. A.},
	title = {Learning variational models with unrolling and bilevel optimization},
	year = {2022},
        eprint={2209.12651},
        archivePrefix={arXiv},
	doi = {10.48550/arxiv.2209.12651},
	note = {arxiv:2209.12651}
}

@article{riccio2024,
title = {Regularization of inverse problems: deep equilibrium models versus bilevel learning},
journal = {Numer. Algebr. Control Optim.},
pages = {},
year = {\noop{3001}in press},
issn = {2155-3289},
doi = {10.3934/naco.2023026},
author = {Danilo Riccio and Matthias J. Ehrhardt and Martin Benning},
}

@book{BerteroBook,
	author = {Bertero, Mario and Boccacci, Patrizia and Ruggiero, Valeria},
	title = {Inverse Imaging with {P}oisson Data},
	publisher = {IOP Publishing},
        address = {Bristol},
	year = {2018},
        doi = {10.1088/2053-2563/aae109},
	isbn = {978-0-7503-1437-4}
}

@article{CombettesWajs2005,
	author = {Combettes, P. L. and Wajs, V. R.},
	title = {Signal Recovery by Proximal Forward-Backward Splitting},
	journal = {Multiscale Model. Simul.},
	volume = {4},
	number = {4},
	pages = {1168--1200},
	year = {2005},
	doi = {10.1137/050626090}
}

@article{ChambollePock2016, 
	author={Chambolle, A. and Pock, T.}, 
	title={An introduction to continuous optimization for imaging}, volume={25}, DOI={10.1017/S096249291600009X}, 
	journal={Acta Numer.}, 
	publisher={Cambridge University Press}, 
	year={2016}, 
	pages={161--319}
}

@article{ChambollePock2011,
	Author = {Chambolle, Antonin and Pock, Thomas},
	Journal = {J. Math. Imaging Vis.},
	Number = {1},
	Pages = {120--145},
	Title = {A First-Order Primal-Dual Algorithm for Convex Problems with Applications to Imaging},
	Volume = {40},
	Year = {2011},
        DOI = {10.1007/s10851-010-0251-1}
}

@article{Boyd2011,
	author = {Boyd, Stephen P. and Parikh, Neal and Chu, Eric and Peleato, Borja and Eckstein, Jonathan},
	journal = {Found. Trends Mach. Learn.},
	number = 1,
	pages = {1--122},
	title = {Distributed Optimization and Statistical Learning via the Alternating Direction Method of Multipliers.},
	volume = 3,
	year = 2011,
        doi = {10.1561/2200000016}
}

@book{Nesterov-2004,
	title = "Introductory lectures on convex optimization: a basic course",
	author = "Nesterov, Y.",
	publisher = "Springer",
	address = "New York",
	year = 2004,
        doi = "10.1007/978-1-4419-8853-9"
}

@ARTICLE{Beck-Teboulle2009,
	AUTHOR =       {A. Beck and M. Teboulle},
	TITLE =        {A fast iterative shrinkage-thresholding algorithm for linear inverse problems},
	JOURNAL =      {SIAM J. Imaging Sci.},
	YEAR =         {2009},
	volume =       {2},
	pages =        {183--202},
        doi = {10.1137/080716542}
}

@article{daubechies2004iterative,
	author = {Daubechies, I. and Defrise, M. and De Mol, C.},
	journal = {Commun. Pure Appl. Math.},
	number = 11,
	pages = {1413--1457},
	publisher = {Wiley Online Library},
	title = {An iterative thresholding algorithm for linear inverse problems with a sparsity constraint},
	volume = 57,
	year = 2004,
        doi = {10.1002/cpa.20042}
}

@ARTICLE{Candes2006,
	author={E. J. {Candes} and J. {Romberg} and T. {Tao}},
	journal={IEEE Trans. Inform. Theory}, 
	title={Robust uncertainty principles: exact signal reconstruction from highly incomplete frequency information},  year={2006},
	volume={52},
	number={2},
        pages={489--509},
        doi={10.1109/TIT.2005.862083}
}

@ARTICLE{Donoho2006,
  author={Donoho, D.L.},
  journal={IEEE Transactions on Information Theory}, 
  title={Compressed sensing}, 
  year={2006},
  volume={52},
  number={4},
  pages={1289-1306},
  doi={10.1109/TIT.2006.871582}}

@article{Stuart_2010, 
	title={Inverse problems: A {B}ayesian perspective}, 
	volume={19}, 
	DOI={10.1017/S0962492910000061}, 
	journal={Acta Numer.}, 
	author={Stuart, A. M.}, 
	year={2010}, 
	pages={451–559}
}

@book{ChanShen,
	author = {Chan, T. and Shen, J.},
	title = {Image Processing and Analysis:  Variational, PDE, Wavelet,
	and Stochastic Methods},
	publisher = {SIAM},
	year = {2005},
	doi = {10.1137/1.9780898717877},
	address =      {Philadelphia},
	edition   = {},
}

@article{Hansen1992,
	author = {Hansen, Per Christian},
	title = {Analysis of Discrete Ill-Posed Problems by Means of the {L}-Curve},
	journal = {SIAM Rev.},
	volume = {34},
	number = {4},
	pages = {561-580},
	year = {1992},
	doi = {10.1137/1034115},
}

@article{Morozov1966,
	title={On the solution of functional equations by the method of regularization},
	author={Vladimir Alekseevich Morozov},
	journal={Dokl. Akad. Nauk SSSR},
	year={1966},
	volume={167},
        number={3},
	pages={510--512}
}

@article{Bertero_2010,
	doi = {10.1088/0266-5611/26/10/105004},
	year = {2010},
	volume = {26},
	number = {10},
	pages = {105004},
	author = {M Bertero and P Boccacci and G Talenti and R Zanella and L Zanni},
	title = {A discrepancy principle for {P}oisson data},
	journal = {Inverse Probl.}
}

@misc{hershey2014deep,
      title={Deep unfolding: Model-based inspiration of novel deep architectures}, 
      author={John R. Hershey and Jonathan Le Roux and Felix Weninger},
      year={2014},
      eprint={1409.2574},
      archivePrefix={arXiv},
      primaryClass={cs.LG},
      note={arXiv:1409.2574},
      doi = {10.48550/arxiv.1409.2574}
}

@article{SIG-111,
	year = {2022},
	volume = {15},
	journal = {Found. Trends Signal Process.},
	title = {Bilevel Methods for Image Reconstruction},
	issn = {1932-8346},
	number = {2-3},
	author = {Caroline Crockett and Jeffrey A. Fessler},
        pages = {121--289},
        doi = {10.1561/2000000111}
}

@article{JCCarolaBilevel2013,
	title = {Image denoising: Learning the noise model via nonsmooth {PDE}-constrained optimization},
	journal = {Inverse Probl. Imaging},
	volume = {7},
	number = {4},
	year = {2013},
	author = {De los Reyes, J. C. and Sch\"{o}nlieb, C.-B.},
        pages = {1183--1214},
        doi = {10.3934/ipi.2013.7.1183}
}

@inbook{CalatroniBilevel2017,
	author = {L. Calatroni and Cao, C. and De los Reyes, J. C. and Sch\"{o}nlieb, C.-B.  and Valkonen, T.},
	title = {Bilevel approaches for learning of variational imaging models},
	year = {2017},
	booktitle = {Variational Methods in Imaging and Geometric Control},
        doi = {10.1515/9783110430394-008},
        pages    = {252--290},
        publisher= {De Gruyter},
        address = {Berlin},
        volume = {18}
}

@article{Kunisch2013,
	author = {Kunisch, K. and Pock, T.},
	title = {A Bilevel Optimization Approach for Parameter Learning in Variational Models},
	journal = {SIAM J. Imaging Sci.},
	volume = {6},
	number = {2},
	year = {2013},
        pages = {938--983},
        doi = {10.1137/120882706}
}

@article{Ochs2016bil,
	author = {Ochs, P. and Ranftl, R.  and Brox, T. and Pock, T.},
	title = {Techniques for Gradient-Based Bilevel Optimization with Non-Smooth Lower Level Problems},
	year = {2016},
	issue_date = {October 2016},
	publisher = {Kluwer Academic Publishers},
	address = {USA},
	volume = {56},
	number = {2},
	issn = {0924-9907},
	doi = {10.1007/s10851-016-0663-7},
	journal = {J. Math. Imaging Vis.},
	pages = {175–194},
	numpages = {20}
}

@InProceedings{pedregosa16,
  title = 	 {Hyperparameter optimization with approximate gradient},
  author = 	 {Pedregosa, Fabian},
  booktitle ={Proceedings of The 33rd International Conference on Machine Learning},
  pages = 	 {737--746},
  year = 	 {2016},
  volume =  {48},
  address =  {New York},
  publisher =    {PMLR}
}

@InProceedings{Grazzi2021,
	title = 	 {Convergence Properties of Stochastic Hypergradients },
	author =       {Grazzi, R. and Pontil, M. and Salzo, S.},
	booktitle = {Proceedings of The 24th International Conference on Artificial Intelligence and Statistics},
	pages = 	 {3826--3834},
	year = 	 {2021},
	volume = 	 {130},
	publisher =    {PMLR},
    address = {New York}
}

@article{arridge2019, 
	title={Solving inverse problems using data-driven models}, 
	volume={28},
	DOI={10.1017/S0962492919000059}, 
	journal={Acta Numer.}, 
	publisher={Cambridge University Press}, 
	author={Arridge, S. and Maass, P. and Öktem, O. and Sch\"{o}nlieb, C.-B.},
	year={2019},
	pages={1–174}}

@article{Sage2019,
	Author = {Sage, Daniel and Pham, Thanh-An and Babcock, Hazen and Lukes, Tomas and Pengo, Thomas and Chao, Jerry and Velmurugan, Ramraj and Herbert, Alex and Agrawal, Anurag and Colabrese, Silvia and Wheeler, Ann and Archetti, Anna and Rieger, Bernd and Ober, Raimund and Hagen, Guy M. and Sibarita, Jean-Baptiste and Ries, Jonas and Henriques, Ricardo and Unser, Michael and Holden, Seamus},
	Journal = {Nat. Methods},
	Number = {5},
	Pages = {387--395},
	Title = {Super-resolution fight club: assessment of {2D} and {3D} single-molecule localization microscopy software},
	Volume = {16},
	Year = {2019},
        DOI = {10.1038/s41592-019-0364-4}
}

@ARTICLE{Lucas2018,
	author={Lucas, A.  and Iliadis, M.  and Molina, R.  and Katsaggelos, A. K.},
	journal={IEEE Signal Process. Mag.}, 
	title={Using Deep Neural Networks for Inverse Problems in Imaging: Beyond Analytical Methods}, 
	year={2018},
	volume={35},
	number={1},
	pages={20-36},
	doi={10.1109/MSP.2017.2760358}}

@article{Nehme_18,
author = {Elias Nehme and Lucien E. Weiss and Tomer Michaeli and Yoav Shechtman},
journal = {Optica},
number = {4},
pages = {458--464},
publisher = {Optica Publishing Group},
title = {Deep-STORM: super-resolution single-molecule microscopy by deep learning},
volume = {5},
year = {2018},
url = {https://opg.optica.org/optica/abstract.cfm?URI=optica-5-4-458},
doi = {10.1364/OPTICA.5.000458}
}

@ARTICLE{Monga2021,
	author={Monga, V. and Li, Y.  and Eldar, Y. C.},
	journal={IEEE Signal Process. Mag.}, 
	title={Algorithm Unrolling: Interpretable, Efficient Deep Learning for Signal and Image Processing}, 
	year={2021},
	volume={38},
	number={2},
	pages={18-44},
	doi={10.1109/MSP.2020.3016905}}

@ARTICLE{Luisier2011,
	author={Luisier, F. and Blu, T. and Unser, M.},
	journal={IEEE Trans. Image Process.}, 
	title={Image Denoising in Mixed {P}oisson-{G}aussian Noise}, 
	year={2011},
	volume={20},
	number={3},
	pages={696--708},
	doi={10.1109/TIP.2010.2073477}}

@article{Soubies_2015_CEL0,
    author = {Soubies, Emmanuel and Blanc-F\'{e}raud, Laure and Aubert, Gilles},
	doi = {10.1137/151003714},
	journal = {SIAM J. Imaging Sci.},
	number = {3},
	pages = {1607--1639},
	title = {A Continuous Exact $\ell_0$ Penalty ({CEL0}) for Least Squares Regularized Problem},
	volume = {8},
	year = {2015}
}

@article{Soubies2017,
author = {Soubies, Emmanuel and Blanc-F\'{e}raud, Laure and Aubert, Gilles},
title = {A Unified View of Exact Continuous Penalties for $\ell\_2$-$\ell\_0$ Minimization},
journal = {SIAM J. Optim.},
volume = {27},
number = {3},
pages = {2034--2060},
year = {2017},
doi = {10.1137/16M1059333}
}

@inproceedings{Stergiopoulou_2023_Col0rmeSSVM,
	author = {Stergiopoulou, Vasiliki and Mukherjee, Subhadip and Calatroni, Luca and Blanc-F{\'e}raud, Laure},
	booktitle = {SSVM 2023: Scale Space and Variational Methods in Computer Vision},
	pages = {498--510},
	publisher = {Springer},
	title = {Fluctuation-Based Deconvolution in Fluorescence Microscopy Using Plug-and-Play Denoisers},
	year = {2023},
        volume = {14009},
        address = {New York},
        doi = {10.1007/978-3-031-31975-4_38}
    }

@incollection{Bertero2008,
	author       = {M. Bertero and H. Lant\'eri and L. Zanni},
	title        = {Iterative image reconstruction: a point of view},
	booktitle    = {Mathematical Methods in Biomedical Imaging and Intensity-Modulated Radiation Therapy ({IMRT})},
	pages        = {37--63},
	publisher    = {Birkhauser-Verlag},
	year         = {2008},
	address      = {Pisa, Italy}
}

\end{document}